\date{November 1999}
\documentclass[a4paper,12pt]{article}
\usepackage{fullpage}
\usepackage{amsmath,amssymb,amscd,amsthm,latexsym}
\usepackage{pictex}

\theoremstyle{plain}  

\newtheorem{thm}{Theorem}[section]
\newtheorem*{thm*}{Theorem}

\newtheorem{lemma}[thm]{Lemma}
\newtheorem{prop}[thm]{Proposition}

\theoremstyle{definition}

\newtheorem{defn}[thm]{Definition}

\newtheorem{rem}[thm]{Remark}

\numberwithin{equation}{section}

\newcommand{\thmref}[1]{Theorem~\ref{#1}}
\newcommand{\propref}[1]{Proposition~\ref{#1}}
\newcommand{\secref}[1]{Section~\ref{#1}}

\newcommand{\lemref}[1]{Lemma~\ref{#1}}

\renewcommand{\leq}{\leqslant}
\renewcommand{\geq}{\geqslant}
\newcommand{\R}{\mathbb{R}}
\newcommand{\Z}{\mathbb{Z}}
\newcommand{\C}{\mathbb{C}}

\newcommand{\calo}{\mathcal{O}}
\newcommand{\dbar}{\bar{\partial}}
\renewcommand{\d}{\mathrm{d}}
\newcommand{\into}{\hookrightarrow}
\renewcommand{\Sp}{\mathrm{Sp}}
\newcommand{\SU}{\mathrm{SU}}
\newcommand{\U}{\mathrm{U}}
\newcommand{\SL}{\mathrm{SL}}

\newcommand{\Or}{\mathrm{O}}
\newcommand{\SO}{\mathrm{SO}}
\newcommand{\abs}[1]{\lvert#1\rvert}
\newcommand{\norm}[1]{\lVert#1\rVert}
\newcommand{\lie}{\mathfrak}
\newcommand{\half}{\tfrac{1}{2}}

\newcommand{\tKO}{\widetilde{KO}}

\newcommand{\vol}{\mathrm{vol}}
\newcommand{\e}{\mathrm{e}}
\newcommand{\ii}{\mathit{i}}
\newcommand{\Id}{\mathrm{Id}}

\DeclareMathOperator{\Jac}{Jac}
\DeclareMathOperator{\ad}{ad}
\DeclareMathOperator{\Ad}{Ad}

\DeclareMathOperator{\tr}{tr}
\DeclareMathOperator{\rk}{rk}

\DeclareMathOperator{\Hom}{Hom}
\DeclareMathOperator{\End}{End}

\begin{document}

\title{Components of spaces of representations and
  stable triples}
\author{Peter B. Gothen}
\date{November 1999}

\maketitle

\begin{abstract}
  We consider the moduli spaces of representations of the fundamental
  group of a surface of genus $g \geq 2$ in the Lie groups $\SU(2,2)$
  and $\Sp(4,\R)$.  It is well known that there is a characteristic
  number, $d$, of such a representation, satisfying the inequality
  $\abs{d} \leq 2g-2$.  This allows one to write the moduli space as a
  union of subspaces indexed by $d$, each of which is a union of
  connected components.  The main result of this paper is that the
  subspaces corresponding to $d = \pm(2g-2)$ are connected in the case
  of representations in $\SU(2,2)$, while they break up into $3\cdot
  2^{2g} + 2g - 4$ connected components in the case of representations in
  $\Sp(4,\R)$.  We obtain our results using the interpretation of the
  moduli space of representations as a moduli space of Higgs bundles,
  and an important step is an identification of certain subspaces
  as moduli spaces of stable triples, as studied by Bradlow and
  Garc{\'\i}a-Prada.
\end{abstract} 

\section{Introduction}

Let $\Sigma$ be a closed Riemann surface of genus $g \geq 2$ and let
$G$ be a connected Lie group.  Consider the space of reductive
representations of the fundamental group of $\Sigma$ in $G$ modulo the
action of $G$ by conjugation,
\begin{displaymath}
  \mathcal{M}_{G} = \Hom(\pi_1(\Sigma),G)^{+}/G,
\end{displaymath}
the superscript ``$+$'' indicating reductive representations.
As is well known, $\mathcal{M}_G$ can also be identified with the moduli
space of reductive flat $G$-bundles over $\Sigma$ and it has an
algebro-geometric interpretation as a moduli space of Higgs bundles
(see Hitchin \cite{hitchin,hitchinduke}).

In this paper we study the connected components of $\mathcal{M}_{G}$
in the cases $G = \Sp(4,\R)$ and $G = \SU(2,2)$.  Previous work on
this type of problem includes the determination of the number of
connected components of $\mathcal{M}_G$ for the groups
$\mathrm{PSL(2,\R)}$ and $\mathrm{PSL(2,\C)}$ by Goldman
\cite{goldman-top}, for the groups $\mathrm{PSL}(n,\R)$ ($n \geq 3$)
by Hitchin \cite{hitchinlie},
 and for the groups $\mathrm{PGL(2,\R)}$, $\mathrm{PU}(2,1)$, and 
 $\mathrm{U}(p,1)$ by
Xia \cite{xia,xia:pu21,xia:up1}.

The first observation is that there is a characteristic number, $d$,
which comes from a characteristic class of the bundle obtained from a
reduction of structure group to the maximal compact subgroups $\U(2)
\subseteq \Sp(4,\R)$ and $\mathrm{S}(\U(2)\times\U(2)) \subseteq
\SU(2,2)$, respectively.  It is well known that this satisfies the
Milnor-Wood type inequality $\abs{d} \leq 2g-2$ (cf.\ 
\secref{sec:milnorwood}).  This allows one to write
\begin{displaymath}
  \mathcal{M}_G = \mathcal{M}_{-(2g-2)} \cup \cdots \cup
  \mathcal{M}_{2g-2},
\end{displaymath}
where each $\mathcal{M}_d$ is a union of connected components.  We can
then state our main result as follows.
\begin{thm*}
  The subspaces $\mathcal{M}_0 \subseteq \mathcal{M}_{G}$ are
  connected for $G = \Sp(4,\R)$ and $G = \SU(2,2)$ and the subspaces
  $\mathcal{M}_{\pm (2g-2)} \subseteq \mathcal{M}_{\SU(2,2)}$ are
  connected.  The subspaces $\mathcal{M}_{\pm (2g-2)} \subseteq
  \mathcal{M}_{\Sp(4,\R)}$ have $3\cdot 2^{2g} + 2g - 4$ connected
  components.
\end{thm*}
The most remarkable aspect of this result is that $\mathcal{M}_{\pm
  (2g-2)} \subseteq \Sp(4,\R)$ breaks up into a number of
different connected components, which are not detected by the first
Chern class given by reduction of structure group to $\U(2) \subseteq
\Sp(4,\R)$.  It seems likely that the remaining $\mathcal{M}_d$ are
connected and we hope to come back to this question on a later occasion.

The method we use for studying the connected components is via the
algebro-geometric interpretation of $\mathcal{M}_G$ as a moduli space
of Higgs bundles, due to Hitchin \cite{hitchin,hitchinlie}.  (A Higgs
bundle is a pair $(E,\Phi)$, where $E$ is a holomorphic rank $n$
degree $d$ vector bundle and $\Phi \in H^0(\Sigma; \End(E)_0 \otimes
K)$, see \secref{sec:higgs} for more details.)  From this point of
view one can define a Hamiltonian circle action on the moduli space
and one uses a moment map for this action as a Morse function, in the
sense of Bott, to obtain topological information about the space (cf.\ 
Hitchin \cite{hitchin,hitchinlie} and Gothen \cite{higgs}).  The
central point is then to identify the critical submanifolds of the
Morse function and to obtain topological information about them.  In
particular, to obtain information about connected components, one
needs to consider the local minima of the Morse function.  It should
be remarked that the moduli spaces have singularities and so one
cannot directly apply Morse theory, however, in the case of the
determination of connected components this difficulty can be
circumvented (see Sections \ref{sec:morse} and
\ref{sec:morseindices}).

In this paper we show that certain critical submanifolds,
corresponding to local minima of the Morse function, can be identified
with moduli spaces of \emph{stable triples}, or spaces closely related
to them, as studied by Bradlow and Garc{\'\i}a-Prada
\cite{bradlow-garcia-prada,garcia-prada}.  In the cases $d = 0$ and
$\abs{d} = 2g-2$ the structure of the moduli spaces of triples is
particularly simple and this allows us to prove the theorem above.  In
the case of $G = \Sp(4,\R)$ and $\abs{d} = 2g-2$ we further need to
use a spectral curve (see Hitchin \cite{hitchinduke}) which is an
unramified covering of $\Sigma$ and the mod 2 index theorem of
Atiyah-Singer to identify the local minima of the Morse function as
certain Prym varieties associated to the covering of $\Sigma$.

This paper is organized as follows.  In \secref{sec:basics} we recall
the basics of the theory of Higgs bundles and their relation to
representations of the fundamental group of a surface in a non-compact
Lie group.  We also recall the concept of a \emph{$Q$-bundle}, of
which a holomorphic triple is a special case, and prove a theorem
(\thmref{thm:quiver}) which is essential for the identification of the
subspace of local minima with a moduli space of holomorphic triples.
Finally we describe the Morse theory on the moduli space and, in
particular, we describe how to find the Morse indices.  It was
observed by Hausel \cite{hausel} that our results on the Morse
indices, together with a theorem of his, imply a theorem of Laumon
\cite{laumon} in this context: the nilpotent cone in the moduli space
$\mathcal{M}$ of rank $n$ Higgs bundles is a Lagrangian subvariety
with respect to the holomorphic symplectic form on $\mathcal{M}$.  We
end this section by briefly describing this.  In
\secref{sec:milnorwood} we reprove the known bound $\abs{d} \leq
2g-2$, using Higgs bundles; we include the proof because it gives some
extra information which is important later on (cf.\ 
\propref{prop:iso}).  In \secref{sec:minima} we analyze the local
minima of the Morse function on the space $\mathcal{M}_G$ for $G =
\Sp(4,\R)$ and $G = \SU(2,2)$ in detail.  Finally, in
\secref{sec:count}, we finish the proof of our main theorem, using the
previous results.

\emph{Acknowledgments.}  Part of this paper is based on my Ph.D.\ 
thesis and I would like to thank my supervisor Nigel Hitchin.
I also benefited from workshops
organized by the VBAC research group under the European algebraic
geometry networks AGE and Europroj, supported by the EU.  This work
was partially supported by Statens Naturvidenskabelige Forskningsr{\aa}d
(Denmark), and by the Funda{\c c}{\~a}o para a Ci{\^e}ncia e a Tecnologia (Portugal)
through the Centro de Matem{\'a}tica da Universidade do Porto and the
project Praxis 2/2.1/MAT/63/94.

\section{Higgs bundles and the topology of moduli spaces}
\label{sec:basics}

\subsection{Higgs bundles}
\label{sec:higgs}

In this section we review some basic facts about Higgs bundles and set
up notation.  For details see Hitchin \cite{hitchin} and Simpson
\cite{simpson:higgs}.

Let $G_{\C}$ be a complex semi-simple Lie group with Lie algebra
$\lie{g}_{\C}$.  Let $G \subset G_{\C}$ be a maximal compact
subgroup with Lie algebra $\lie{g}$.  Thus there is a compact real
structure $\tau \colon \lie{g}_{\C} \to \lie{g}_{\C}$ whose fixed
point set is $\lie{g}$.  Denoting the $-1$-eigenspace of $\tau$ by
$\lie{g}^{\perp}$ we then have $\lie{g}_{\C} = \lie{g} \oplus
\lie{g}^{\perp}$. 

Non-abelian Hodge
theory gives an equivalence between reductive representations of
$\pi_1(\Sigma)$ in $G_{\C}$ and Higgs bundles over $\Sigma$, which we
now describe.  Let 
\begin{displaymath}
  \rho \colon \pi_1(\Sigma) \to G_{\C}
\end{displaymath}
be a reductive representation.  This data is equivalent to having a
principal bundle 
\begin{displaymath}
  P_{\C} = \tilde{\Sigma}\times_{\rho}G_{\C}
\end{displaymath}
with a reductive flat connection $B \in \Omega^1(P_{\C};\lie{g}_{\C})$
(here $\tilde{\Sigma}$ is the universal cover of $\Sigma$).

If we have a
metric in $P_{\C}$, i.e.\ a reduction of structure group from $G_{\C}$
to $G$, we can write 
\begin{displaymath}
  i^*B = A + \theta
\end{displaymath}
where $i \colon P \into P_{\C}$ is the inclusion of the principal
$G$-bundle $P$ given by the reduction of structure group, $A$ is a
connection on $P$, and $\theta \in \Omega^1(P;\lie{g}^{\perp})$ is a
tensorial form, which can therefore be thought of as an element of
$\Omega^1(\Sigma; P\times_{\Ad} \lie{g}^{\perp})$.  

Given a complex
representation of $G$ (e.g.\ the adjoint representation on
$\lie{g}_{\C}$), we
have the usual decomposition of the covariant derivative $\d_A$ in its
$(1,0)$- and $(0,1)$-parts:
\begin{displaymath}
  \d_A = \partial_A + \dbar_A.
\end{displaymath}
Similarly, we can write 
\begin{displaymath}
  \theta = \Phi - \tau(\Phi),
\end{displaymath}
for a unique $\Phi \in \Omega^{1,0}(\Sigma;\Ad P_{\C})$ (by abuse of
notation, we denote by $\tau$
the combination of the compact real structure $\tau$ on $\lie{g}_{\C}$
and conjugation on the form component).

Corlette \cite{corlette} and Donaldson \cite{Donaldson} proved that
there exists a harmonic metric in $P$, that is, a metric such that
$(A,\Phi)$ obtained via the above procedure satisfy Hitchin's
equations
\begin{align*}
  F(A) - [\Phi,\tau(\Phi)] &= 0, \\
  \dbar_A\Phi &= 0.
\end{align*}
This, in turn, gives a principal \emph{Higgs bundle}, i.e.\ a pair
$(P_{\C},\Phi)$ consisting of a holomorphic principal bundle $P$ (with
holomorphic structure defined by $\dbar_A$) and a \emph{Higgs field}
$\Phi \in H^0(\Sigma;\Ad P_{\C}\otimes K)$, where $K$ denotes the
canonical bundle of $\Sigma$.  
Given a representation, $V$, of
$G_{\C}$ one then obtains a Higgs vector bundle $(E,\Phi)$, where $E =
P_{\C} \times_{G_{\C}} V$ and $\Phi \in H^0(\Sigma;\End(E)\otimes
K)$.  The two main examples we have in mind are the adjoint
representation $V = \lie{g}_{\C}$ and the fundamental representation
of $G_{\C} = \SL(n,\C)$.
If the original representation $\rho$ of $\pi_1(\Sigma)$ is
irreducible then $(E,\Phi)$ is \emph{stable}, i.e.\
\begin{displaymath}
  \mu(F) < \mu(E)
\end{displaymath}
for any proper non-trivial $\Phi$-invariant subbundle $F$ of $E$ (here
$\mu(E) = \deg(E)/\rk(E)$ is the \emph{slope} of the holomorphic
bundle $E$).  Allowing equality in the above inequality gives the
notion of a semi-stable Higgs bundle.  Finally, if $\rho$ is
reductive, then the corresponding Higgs bundle is poly-stable, i.e.\ 
it is a direct sum of lower rank Higgs bundles, all of the same slope.

By a theorem of Hitchin \cite{hitchin} and Simpson \cite{simpson},
the above procedure can be reversed, by finding a harmonic metric in
the Higgs bundle.  This produces a reductive
representation of $\pi_1(\Sigma)$ from a poly-stable Higgs bundle.
This gives a homeomorphism
\begin{displaymath}
  \mathcal{M}_{G_{\C}} \to
  \Hom(\pi_1(\Sigma),G_{\C})^{+}/G_{\C}.
\end{displaymath}
where $\mathcal{M}_{G_{\C}}$ is the moduli space of poly-stable
principal $G_{\C}$ Higgs bundles.

We finish by recalling the description of the Zariski tangent space to
$\mathcal{M}_{G_{\C}}$ at $(P_{\C},\Phi)$ given by Biswas and Ramanan
\cite{biswasramanan}.  This is the first hyper-cohomology
$\mathbb{H}^1(C_{\C}^{\bullet})$ of the
complex of sheaves
\begin{equation}
  \label{eq:cxtangent}
  C_{\C}^{\bullet} : \mathcal{O}(\Ad P_{\C}) \xrightarrow{\ad(\Phi)}
  \mathcal{O}(\Ad P \otimes K).
\end{equation}
From this they deduce the long exact sequence
\begin{multline}
  \label{eq:cxlong}
  0 \to \mathbb{H}^0(C_{\C}^{\bullet}) \to H^0(\Sigma;\Ad P_{\C}) \to
  H^0(\Sigma;\Ad P \otimes K)
  \to T_{(P_{\C},\Phi)}\mathcal{M}_{G_{\C}} \\
  \to H^1(\Sigma;\Ad P_{\C}) \to H^1(\Sigma;\Ad P \otimes K) \to
  \mathbb{H}^2(C_{\C}^{\bullet}) \to 0.
\end{multline}
$(P_{\C},\Phi)$ is a smooth point of the moduli space if
$\mathbb{H}^0(C_{\C}^{\bullet})$ and $\mathbb{H}^2(C_{\C}^{\bullet})$
vanish (by Serre duality, it is sufficient to check that
$\mathbb{H}^0(C_{\C}^{\bullet}) = 0$).  From this one sees that stable
Higgs bundles represent smooth points of the moduli space.  The
dimension of $\mathcal{M}_{G_{\C}}$ can be calculated using
Riemann-Roch to be
\begin{displaymath}
  \dim_{\C}\mathcal{M}_{G_{\C}}
  = \chi(\mathcal{O}(\Ad P \otimes K))
  - \chi(\mathcal{O}(\Ad P_{\C}))
  = 2\dim_{\C} \lie{g}_{\C}(g-1)
\end{displaymath}

\subsection{Real groups}

Hitchin \cite{hitchin,hitchinlie} showed how to use Higgs bundles to
study representations of $\pi_1(\Sigma)$ in real (non-compact) Lie
groups.  Next we recall the relevant parts of this theory.

Let $G_r \subset G_{\C}$ be a real form, given by a real structure
$\sigma \colon \lie{g}_{\C} \to \lie{g}_{\C}$.  Let $K \subset G_{r}$
be a maximal compact subgroup and let $\lie{k} \subset \lie{g}_{r}$ be
the corresponding inclusion of Lie algebras.  Let $\lie{k}^{\perp}$ be
the orthogonal complement to $\lie{k}$ with respect to the Killing
form, then we can write $\lie{g}_{r} = \lie{k} \oplus
\lie{k}^{\perp}$, where the Killing form is negative definite on
$\lie{k}$ and positive definite on $\lie{k}^{\perp}$.  Define a
complex linear involution $\phi \colon \lie{g}_{\C} \to \lie{g}_{\C}$
by $\phi_{\vert \lie{k}_{\C}} = 1$ and $\phi_{\vert
\lie{k}^{\perp}_{\C}} = -1$.  Define another real structure on
$\lie{g}_{\C}$ by $\tau = \sigma \phi = \phi \sigma$.  It is then 
easy to see the that the corresponding real subgroup $G \subset 
G_{\C}$ is a maximal compact subgroup and, clearly, $\lie{k} = \lie{g} 
\cap \lie{g}_r$.  

Now suppose that we have a reductive representation of $\pi_1(\Sigma)$
in $G_r$ and let $B$ be the associated flat connection on the
principal $G_r$-bundle $P_{G_r}$.  The theorem of Donaldson and
Corlette also applies in this case and gives a reduction of structure
group to $K \subset G_r$: let $i \colon P_K \into P_{G_{\C}}$ be the
inclusion of principal bundles given by combining the reduction of
structure group with the inclusion $G_r \subset G_{\C}$.  In the
decomposition $i^*B = A + \theta$, $A$ and $\theta$ will be fixed by
$\sigma$, while $\tau(A) = A$ and $\tau(\theta) = -\theta$.  Thus
$\dbar_A$ is fixed by $\phi$, and $\Phi$ is in the $-1$-eigenspace of
$\phi$.  This means that the corresponding Higgs bundle is of
the form 
\begin{math}
  (P_{K_{\C}},\Phi),
\end{math}
where 
\begin{itemize}
\item $P_{K_{\C}}$ is a holomorphic principal $K_{\C}$-bundle,  
\item $\Phi \in H^0(\Sigma; \Ad_{\lie{k}_{\C}^{\perp}} P_{K}\otimes
  K)$ (where we use the notation $\Ad_{\lie{k}_{\C}^{\perp}} P_{K} =
  P_K \times_{\Ad} \lie{k}_{\C}^{\perp}$).
\end{itemize}
Conversely, such a Higgs bundle gives a representation of
$\pi_1(\Sigma)$ in $G_r$.  We then have a homeomorphism
\begin{displaymath}
  \mathcal{M}_{G_{r}} \to
  \Hom(\pi_1(\Sigma),G_{r})^{+}/G_{r},
\end{displaymath}
where $\mathcal{M}_{G_{r}}$ is the moduli space of poly-stable Higgs
bundles of the above type.  Alternatively $\mathcal{M}_{G_r}$ can be
thought of as the moduli space of solutions $(A,\Phi)$ to Hitchin's
equations modulo $K$-gauge equivalence: then $A$ is a connection on a
principal $K$-bundle $P_K$ and $\Phi \in 
\Omega^{1,0}(\Sigma;\Ad_{\lie{k}_{\C}^{\perp}} P_{K})$.

The analogue to \eqref{eq:cxtangent} in this context is that the
Zariski tangent space to $\mathcal{M}_{G_{r}}$ is the first
hyper-cohomology of the complex of sheaves
\begin{equation}
  \label{eq:realtangent}
  C_r^{\bullet} :
  \mathcal{O}(\Ad_{\lie{k}_{\C}} P_{K})
  \xrightarrow{\ad(\Phi)}
  \mathcal{O}(\Ad_{\lie{k}_{\C}^{\perp}} P_{K}\otimes K),
\end{equation}
where we use the notation $\Ad_{\lie{k}_{\C}} P_{K} = P_K \times_{\Ad}
\lie{k}_{\C} = \Ad P_{K_{\C}}$.
The analogue to the long exact sequence \eqref{eq:cxlong} is
\begin{multline}
  \label{eq:reallong}
  0 \to \mathbb{H}^0(C_{r}^{\bullet}) \to
  H^0(\Sigma;\Ad_{\lie{k}_{\C}} P_{K}) \to
  H^0(\Sigma;\Ad_{\lie{k}_{\C}^{\perp}} P_{K} \otimes K)
  \to T_{(P_{K},\Phi)}\mathcal{M}_{G_{r}} \\
  \to H^1(\Sigma;\Ad_{\lie{k}_{\C}} P_{K}) \to
  H^1(\Sigma;\Ad_{\lie{k}_{\C}^{\perp}} P_{K} \otimes K) \to
  \mathbb{H}^2(C_{r}^{\bullet}) \to 0.
\end{multline}
The smooth points of the moduli space
are those for which $\mathbb{H}^0(C_r^{\bullet}) =
\mathbb{H}^2(C_r^{\bullet}) = 0$ and again the stable Higgs bundles
represent smooth points.
The dimension of $\mathcal{M}_{G_r}$ can be calculated as before to be
\begin{displaymath}
  \dim_{\C}\mathcal{M}_{G_r} = \dim_{\C} \lie{g}_{\C}(g-1) = \half
  \dim_{\C}\mathcal{M}_{G_{\C}}.
\end{displaymath}
We finish this section by giving two examples of this setup.  First
consider $G_r = \SU(n,n)$ which is a real form of $\SL(2n,\C)$.
The Higgs vector bundles $(E,\Phi)$ corresponding to representations of
$\pi_1(\Sigma)$ in $\SU(n,n)$ are of the form
\begin{equation}\label{eq:suhiggsbundles}
        E = V \oplus V'
        \quad \text{and} \quad 
        \Phi = 
        \begin{pmatrix}
                0       &       b       \\
                c       &       0
        \end{pmatrix},  
\end{equation}
where $V$ and $V'$ are rank $n$ vector bundles with $\Lambda^n V
\otimes \Lambda^nV' \cong \mathcal{O}$, $b \in H^0(\Hom(V',V) \otimes
K)$, and $c \in H^0(\Hom(V,V') \otimes K)$.  Two $\SU(n,n)$
representations are conjugate if and only if
the corresponding Higgs bundles of this form are isomorphic by an
isomorphism which is of the form 
\begin{math}
  \left(
    \begin{smallmatrix}
      g & 0 \\
      0 & g'
    \end{smallmatrix}
  \right)
\end{math}
and of determinant one.

The second example is $G_r = \Sp(2n,\R)$;
this is a split real form of $\Sp(2n,\C)$.  The Higgs vector bundles
$(E,\Phi)$ obtained from the standard representation of $\Sp(2n,\C)$ on
$\C^{2n}$, and corresponding to representations of $\pi_1(\Sigma)$ in
$\Sp(2n,\R)$ are of the form
\begin{equation}\label{eq:sphiggsbundles}
        E = V \oplus V^*
        \quad \text{and} \quad 
        \Phi = 
        \begin{pmatrix}
                0       &       b       \\
                c       &       0
        \end{pmatrix},  
\end{equation}
where $V$ is a rank $n$ vector bundle, $b \in H^0(S^2 V \otimes K)$,
and $c \in H^0(S^2 V^* \otimes K)$.  Two $\Sp(2n,\R)$
representations are conjugate if and only if
the corresponding Higgs bundles of this form are isomorphic by an
isomorphism which is of the form 
\begin{math}
  \left(
    \begin{smallmatrix}
      g & 0 \\
      0 & g^{t}
    \end{smallmatrix}
  \right)
\end{math}.

\subsection{$Q$-bundles and triples}
\label{sec:q-bundles-triples}

The special forms \eqref{eq:suhiggsbundles} and
\eqref{eq:sphiggsbundles} suggest a different point of view, that of
$Q$-bundles.  This notion, due to Alastair King,
provides a general framework for considering a large number
of the various kinds of vector bundles with extra structure, which
have been studied in recent years.  The vortex pairs of Bradlow
\cite{bradlow}, the triples of Garc{\'\i}a-Prada, introduced in
\cite{garcia-prada} and studied systematically by him and Bradlow in 
\cite{bradlow-garcia-prada} (and also Higgs bundles), are all examples
of $Q$-bundles.

Let $Q$ be a quiver, that is, $Q$ is a directed graph, specified by a
set of vertices $Q_0$ and a set of arrows $Q_1$, together with head
and tail maps $h,t \colon Q_1 \to Q_0$.

\begin{defn}
A \emph{$Q$-bundle\/} over a Riemann surface $\Sigma$ is a collection of
holomorphic vector bundles $\{E_i\}_{i \in Q_0}$ over $\Sigma$ and a
collection of holomorphic maps 
$\{\phi_a \colon E_{t(a)} \to E_{h(a)} \}_{a \in Q_1}$.
A \emph{twisted $Q$-bundle} is given by in addition specifying a
linebundle $L_a$ for each arrow $a$.  The maps $\phi_a$ are then
required to be holomorphic maps
\begin{math}
\phi_a \colon E_{t(a)} \to E_{h(a)} \otimes L_a
\end{math}.
\end{defn}

\begin{figure}[h]
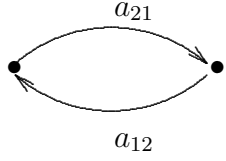

    \centering
\font\thinlinefont=cmr5
\begingroup\makeatletter\ifx\SetFigFont\undefined%
\gdef\SetFigFont#1#2#3#4#5{%
  \reset@font\fontsize{#1}{#2pt}%
  \fontfamily{#3}\fontseries{#4}\fontshape{#5}%
  \selectfont}%
\fi\endgroup%
\mbox{\beginpicture
\setcoordinatesystem units <1.00000cm,1.00000cm>
\unitlength=1.00000cm
\linethickness=1pt
\setplotsymbol ({\makebox(0,0)[l]{\tencirc\symbol{'160}}})
\setshadesymbol ({\thinlinefont .})
\setlinear
%
%
\linethickness= 0.500pt
\setplotsymbol ({\thinlinefont .})
%
%
\plot  4.847 24.408  5.080 24.289  4.930 24.504 /
\circulararc 82.224 degrees from  5.080 24.289 center at  3.810 22.834
%
%
\linethickness= 0.500pt
\setplotsymbol ({\thinlinefont .})
%
%
\plot  2.773 24.011  2.540 24.130  2.690 23.915 /
\circulararc 82.224 degrees from  2.540 24.130 center at  3.810 25.585
%
%
\put{$\bullet$} [lB] at  2.381 24.130
%
%
\put{$\bullet$} [lB] at  5.080 24.130
%
%
\put{$a_{21}$} [lB] at  3.810 24.924
%
%
\put{$a_{12}$} [lB] at  3.810 23.178
\linethickness=0pt
\putrectangle corners at  2.381 25.076 and  5.097 23.139
\endpicture}
    \caption{The quiver $Q$}
    \label{fig:1}
\end{figure}

We shall only consider $Q$-bundles of a particularly simple form: we
let $Q$ be a quiver with $2$ vertices and exactly one arrow connecting
the vertices in each direction (see fig.\ \ref{fig:1}).

We denote the the arrows by $a_{ij}$, where $a_{ij}$ is the arrow
going from $j$ to $i$.  Also, the maps will be twisted by the
canonical bundle $K$. Thus, from now on, a $Q$-bundle is a pair
\begin{displaymath}
        \mathbf{E}=(\underline{E},\underline{\Phi}),
\end{displaymath}
where $\underline{E}=\{E_1,E_2\}$ and $\underline{\Phi}=\{\phi_{ij}\}$.
Here, each $E_i$ is a holomorphic vector bundle on $\Sigma$ and
$\phi_{ij}$ is a holomorphic section of $\Hom(E_j,E_i \otimes K)$.

A particularly interesting special case occurs when $\phi_{12} = 0$.
The data of the above type of $Q$-bundle then comes down to a triple
$(E_1,E_2,\phi)$, where $\phi \in H^0(\Sigma;\Hom(E_1,E_2)\otimes
K)$.  If we define $\tilde{E}_2 = E_2 \otimes K$ then this is
equivalent to a holomorphic triple $(E_1,\tilde{E}_2,\phi)$ in the
sense of Bradlow and Garc{\'\i}a-Prada \cite{bradlow-garcia-prada}.

Given a $Q$-bundle $\mathbf{E}=(\underline{E},\underline{\Phi})$, we can
define an associated Higgs bundle $(E,\Phi)$ by putting
\begin{equation}
        E = E_1 \oplus E_2
        \quad\text{and}\quad
        \Phi = (\phi_{ij}),
\end{equation}
where $(\phi_{ij})$ is the matrix of $\Phi$ with respect to the above
direct sum decomposition of $E$.  Note that the Higgs bundles of the
form \eqref{eq:suhiggsbundles} or \eqref{eq:sphiggsbundles} arise in
this way.  Conversely, given a Higgs bundle of the special form
\eqref{eq:suhiggsbundles} or \eqref{eq:sphiggsbundles} we get an
associated $Q$-bundle.

There are equations for preferred special metrics in a $Q$-bundle, the
\emph{$Q$-vortex equations}.  Choose a metric compatible with the
complex structure on $\Sigma$ and, for convenience, normalize it so
that $\vol(\Sigma) = 2\pi$.  This of course also gives a Hermitian
metric in the canonical bundle $K$.  The $Q$-vortex equations are
equations for Hermitian metrics in $E_1$ and $E_2$ and in our case
they take the form
\begin{equation}
  \label{eq:qvortex} \left\{
  \begin{split}
    \ii \Lambda F(A_1) + \phi_{12}\phi_{12}^* - \phi_{21}^*\phi_{21}
    &= \tau_1 \Id_{E_1} \\
    \ii \Lambda F(A_2) + \phi_{21}\phi_{21}^* - \phi_{12}^*\phi_{12}
    &= \tau_2 \Id_{E_2} \\
  \end{split}\right.
\end{equation}
where $F(A_i)$ is the curvature of the metric connection in $E_i$,
$\Lambda$ denotes contraction with the K{\"a}hler form of $\Sigma$, and
$\phi_{ij}^*$ denotes the adjoint taken with respect to the metric
obtained from the metrics on $E_i$ and $K$.
The parameters $(\tau_1,\tau_2)$ are real, subject to the condition
\begin{displaymath}
  \sum_{i=1}^2 \bigl(\deg(E_i) - \tau_i \rk(E_i) \bigr) = 0,
\end{displaymath}
obtained by taking traces in the equations \eqref{eq:qvortex}, summing
and integrating over $\Sigma$ (thus there is really only one real
parameter involved, which is usually taken to be $\tau=\tau_2$).
There is a stability condition for $Q$-bundles, such that any
$Q$-bundle which supports a solution to the $Q$-vortex equations is a
direct sum of stable $Q$-bundles.  In our case the condition is
\begin{equation}
  \label{eq:taustability}
  \sum_{i=1}^2 \bigl(\deg(F_i) - \tau_i \rk(F_i) \bigr) < 0
\end{equation}
for any proper $Q$-subbundle $\mathbf{F}$ of $\mathbf{E}$.  Note that
the condition depends on the parameters $(\tau_1,\tau_2)$.  Bradlow
and Garc{\'\i}a-Prada \cite{bradlow-garcia-prada} constructed moduli spaces
of stable triples, varying with the parameter $\tau$.

We shall only need to consider the case $\tau_1 = \tau_2 = \mu(E)$ so
we shall assume this from now on.  The stability condition
\eqref{eq:taustability} can then be reformulated as
\begin{equation}
  \label{eq:alphastability}
  \mu(F) < \mu(E)
\end{equation}
for any proper $Q$-subbundle
$\mathbf{F}=(\{F_1,F_2\},\{\phi_{12},\phi_{21}\})$ of $\mathbf{E}$,
and where we write $F = F_1 \oplus F_2$.

But, obviously, $F = F_1 \oplus F_2 \subset E$ is a $\Phi$-invariant
subbundle, thus stability of the Higgs bundle $(E,\Phi)$ implies
stability of the $Q$-bundle $\mathbf{E}$.  The following lemma will
allow us to conclude that the converse also holds.

\begin{lemma}
  \label{lemma:quiver}
  Let $(E,\Phi)$ be a Higgs bundle of the form $E = E_1 \oplus E_2$ and
  \begin{displaymath}
  \Phi =
    \begin{pmatrix}
      0 & \phi_{12} \\
      \phi_{21} & 0
    \end{pmatrix}.
  \end{displaymath}
  Let $\mathbf{E} = (\{E_1,E_2\} , \{\phi_{12},\phi_{21}\})$ be the
  associated $Q$-bundle.  Let $F' \subset E$ be a $\Phi$-invariant
  subbundle.  Then there is a $Q$-subbundle $\mathbf{E}'=
  (\{E'_1,E'_2\} , \{\phi_{12},\phi_{21}\})$
  of $\mathbf{E}$ such that
  \begin{displaymath}
    \mu(F') \leq \mu(E'),
  \end{displaymath}
  where $E'=E'_1 \oplus E'_2$.
\end{lemma}

\begin{proof}
Let 
$\pi_{i}\colon E \to E_{i}$ be the projection on the $i$th factor.  Let 
$F_{i} \subset E_i$ and $G_{i} \subset F'$ be the subbundles which are
generated by the image and kernel of $\pi_i$, respectively.
Then
$F_1$ and $G_2$ are contained in $E_1$, $F_2$ and $G_1$ are contained in 
$E_2$, and we have sequences of 
vector bundles
\begin{displaymath}
        0 \to G_i \to F' \to F_i \to 0,
\end{displaymath}
which are generically short exact. Hence, 
\begin{math}
        \deg(F') \leq \deg(G_i) +\deg(F_i),
\end{math}
and putting $F = F_1 \oplus F_2$ and $G = G_2 \oplus G_{1}$, it
follows that
\begin{displaymath}
        2 \deg(F') \leq \deg(F) + \deg(G).
\end{displaymath}
Clearly
\begin{math}
        2 \rk(F') = \rk(F) + \rk(G),
\end{math}
so that 
\begin{equation}
  \label{eq:convex}
    \mu(F') \leq \frac{\rk(F)}{\rk(F) + \rk(G)} \mu(F)
    + \frac{\rk(G)}{\rk(F) + \rk(G)} \mu(G),
\end{equation}
and therefore either $\mu(F) \geq \mu(F')$ or $\mu(G) \geq \mu(F')$.
Provided that $F$ and $G$ give $Q$-subbundles of $\mathbf{E}$ we can
then take $\mathbf{E}'$ to be the $Q$-bundle associated to either $F$
or $G$.

It thus remains to see show that $F$ and $G$ are $\Phi$-invariant and,
therefore, define $Q$-sub\-bundles of $\mathbf{E}$.  First, let $x_1
\in F_1$. If we write $x_1 = \pi_{1}(x)$ for some $x = x_1 + x_2$ in
$F'$, then
\begin{displaymath}
        \Phi(x) = \Phi(x_1) + \Phi(x_2).
\end{displaymath}
By our assumption on the matrix for $\Phi$, it follows that $\Phi(x_1)
\in E_2$ and $\Phi(x_2) \in E_1$.  Then $\pi_{1} (\Phi (x)) = \Phi
(x_2) \in E_{1}$ and $\pi_2 (\Phi(x)) = \Phi (x_1) \in E_2$. But
$\Phi(x) \in F'$ because $F'$ is $\Phi$-invariant, and thus $\Phi(x_2)
\in F_{1}$ and $\Phi(x_1) \in F_{2}$. Of course, we can repeat the
argument with $x_2 \in F_2$ and hence, $F$ is $\Phi$-invariant.  The
proof that $G$ is $\Phi$-invariant is similar.  Let $x_1 \in G_2$.  By
assumption, $\Phi(x_1) \in E_{2}$.  But $G_2 \subset F'$, so
$\Phi(x_1) \in F'$ as well.  It follows that $\Phi(x_1) \in G_{1}$ and
thus, $G$ is $\Phi$-invariant.  We have thus seen that $F$ and $G$
define $Q$-subbundles of $\mathbf{E}$ and this finishes the proof.
\end{proof}

As an immediate consequence we have the following theorem.

\begin{thm}
\label{thm:quiver}
Let $Q$ be a quiver with two vertices and one arrow connecting the
vertices in each direction, and let $\mathbf{E} = (\{E_1,E_2\} , \{
\phi_{12},\phi_{21}\})$ be a $Q$-bundle.  Let $(E,\Phi)$ be the
associated Higgs bundle as above; thus $E = E_1 \oplus E_2$ and
\begin{displaymath}
\Phi =
  \begin{pmatrix}
    0 & \phi_{12} \\
    \phi_{21} & 0
  \end{pmatrix}.
\end{displaymath}
Then $\mathbf{E}$ is stable if and only if $(E,\Phi)$ is.
Furthermore, if $(E,\Phi)$ is poly-stable, i.e.\ the direct sum of
lower rank stable Higgs bundles, these lower rank Higgs bundles
are $Q$-subbundles of $\mathbf{E}$.
\end{thm}

\begin{proof}
  The only assertion that requires proof is the final one.  Suppose
  that $(E,\Phi)$ is poly-stable and that $(F',\Phi')$ is a proper
  stable Higgs subbundle of $(E,\Phi)$ with $\mu(F') = \mu(E)$.  By
  semi-stability of $(E,\Phi)$ the bundle called $G$ in the proof of
  \lemref{lemma:quiver} must then satisfy $\mu(G) = \mu(E) = \mu(F')$
  and, since $G \subseteq F'$, it follows by stability of $F'$ that
  $G=F'$.  But $G$ is a $Q$-subbundle so this finishes the proof.
\end{proof}

\subsection{Connected components and Morse theory}
\label{sec:morse}

We shall use Hitchin's method \cite{hitchin,hitchinlie}, which we
shall now review, for finding the connected components of
$\mathcal{M}_{G_{r}}$.  The idea is to use discrete invariants of
flat bundles for dividing $\mathcal{M}_{G_{r}}$ into subspaces which
are unions of connected components.  We then show that these subspaces
are, in fact, connected.  For this consider the moduli space
$\mathcal{M}_{G_r}$ of Higgs bundles as the space of solutions
$(A,\Phi)$ to Hitchin's equations modulo gauge equivalence.  The
function
\begin{align*}
  f \colon \mathcal{M} &\to \R \\
  (A,\Phi) &\mapsto \norm{\Phi}^2 = \int_{\Sigma}\abs{\Phi}^2\d\vol
\end{align*}
is proper.  Thus, a subspace $N$ of $\mathcal{M}$ is connected if the
subspace of local minima of $f$ on $N$ is connected.  

Restrict for a moment attention to irreducible solutions to Hitchin's
equations (i.e.\ stable Higgs bundles); these are smooth points of
$\mathcal{M}$.  In order to identify the subspaces of local minima of
$f$ one uses the fact that it is a moment map for the $S^1$ action on
$\mathcal{M}$, given by $(A,\Phi) \mapsto (A,\e^{\ii\theta}\Phi)$:
this implies that the critical points of $f$ are exactly the fixed
points of the circle action.  Now, $(A,\Phi)$ represents a fixed point
if and only if there is an infinitesimal gauge transformation $\psi
\in \Omega^0(\Sigma; P_K\times_{\Ad}\lie{k})$ such that
\begin{align}
  d_A\psi &=0, \\
  [\psi,\Phi] &= \ii\Phi. \label{eq:psiPhi}
\end{align}

Let $(E,\Phi)$ be a Higgs vector bundle obtained from a complex
representation of $K$, then this can be decomposed in eigenspaces for
the covariantly constant gauge transformation $\psi$.  Thus
\begin{equation}\label{eq:decomp}
  E = \bigoplus_{m} U_m,
\end{equation}
where $\psi_{\vert U_m} = \ii m$.  Then \eqref{eq:psiPhi} shows that
$\Phi \colon U_m \to U_{m+1} \otimes K$.  The case 
\begin{displaymath}
  E = \Ad P_{\C} = \Ad_{\lie{k}_{\C}} P_{K} \oplus
  \Ad_{\lie{k}_{\C}^{\perp}} P_{K},
\end{displaymath}
is of particular interest.  Since the adjoint action of
$\psi \in \lie{k}$ and the involution $\phi$ on $\lie{g}_{\C}$
commute, this decomposition of $\Ad P_{\C}$ is compatible with the
decomposition $\Ad P_{\C} = \bigoplus_{m} U_m$.  From $\ad(\Phi)
\colon U_m \to U_{m+1} \otimes K$ we conclude that $\Phi \in
H^0(\Sigma; U_1 \otimes K)$ and since $\Phi \in
H^0(\Sigma;\Ad_{\lie{k}_{\C}^{\perp}} P_{K})$ it follows that $U_1
\subset \Ad_{\lie{k}_{\C}^{\perp}} P_{K}$.  Furthermore, $\ad(\Phi)$
interchanges $\Ad_{\lie{k}_{\C}} P_{K}$ and
$\Ad_{\lie{k}_{\C}^{\perp}} P_{K}$ and we therefore have
\begin{align*}
  \Ad_{\lie{k}_{\C}} P_{K} &= \bigoplus_{k} U_{2k}, \\
  \Ad_{\lie{k}_{\C}^{\perp}} P_{K}
  &= \bigoplus_{k} U_{2k+1},
\end{align*}
where $k$ is integer.

Two additional pieces of information will be useful.  The first is
that there is an
isomorphism $\Ad P_{\C} \xrightarrow{\cong} \Ad P_{\C}^*$ from the
adjoint bundle to the co-adjoint bundle given by the Killing form on
$\lie{g}_{\C}$ and it is trivial to check that under this
\begin{equation}
  \label{eq:UmtoU-m*}
  U_m \xrightarrow{\cong} U_{-m}^*.
\end{equation}

The second useful piece of information is that one can calculate the
value of the Morse function $f$ at a Higgs bundle of the form
\eqref{eq:decomp}: denoting the component of $\Phi$ mapping $U_i$ to
$U_{i+1} \otimes K$ by $\phi_i$ one shows easily, using Hitchin's
equations, that
\begin{equation}
  \label{eq:vanishing}
  \norm{\phi_i}^2 = \norm{\phi_{i-1}}^2 + \bigl(\deg(U_i) -
  \mu(E)\rk(U_i) \bigr).
\end{equation}
In particular, if $\mu(U_i) = \mu(E)$ for all $i$ then $\Phi$ must
vanish.

Finally, consider the case of local minima of $f$ which are not
represented by stable Higgs bundles.  For simplicity we restrict
attention to bundles of the form \eqref{eq:suhiggsbundles} or
\eqref{eq:sphiggsbundles}.  These are then direct sums of stable Higgs
bundles of lower rank.  But from \thmref{thm:quiver} we can conclude
that these lower rank Higgs bundles decompose as direct sums of
subbundles of the bundles appearing in the direct sum decomposition
of the original Higgs bundle.

\subsection{Morse indices}
\label{sec:morseindices}

Given a poly-stable Higgs bundle $(P_{K_{\C}},\Phi)$ which represents
a critical point of $f$ it is necessary to decide whether it is a
local minimum.  This can be done using the observation of Hitchin
\cite{hitchinlie} that if $\psi$ acts with weight $m$ on an element of
$\Ad_{\lie{k}_{\C}}P_K$ then the corresponding eigenvalue of the
Hessian of $f$ is $-m$, while a weight $m$ on
$(\Ad_{\lie{k}_{\C}}P_K)\otimes K$ gives the eigenvalue $1-m$.  It
follows that the subspace of $T_{(P_{K_{\C}},\Phi)}\mathcal{M}_{G_r}$
on which the Hessian of $f$ is negative definite is
$\mathbb{H}^1(C^{\bullet}_{-})$, where $C^{\bullet}_{-}$ is the
complex of sheaves
\begin{displaymath}
  C^{\bullet}_{-} : \mathcal{O}\Bigl(\bigoplus_{k \geq 1} U_{2k}\Bigr)
  \xrightarrow{\ad(\Phi)} \mathcal{O}
  \Bigl(\bigoplus_{k \geq 1} U_{2k+1} \otimes K\Bigr).
\end{displaymath}
In other words, the Morse index is
$\dim_{\R}\mathbb{H}^1(C^{\bullet}_{-})$.  At a smooth point of
$\mathcal{M}_{G_r}$, where
$\mathbb{H}^0(C^{\bullet}_{-})=\mathbb{H}^2(C^{\bullet}_{-})=0$, the
Morse index can be calculated using the Riemann-Roch theorem:
\begin{multline}\label{eq:morseindex}
  \dim_{\C} \mathbb{H}^1(C^{\bullet}_{-})
  = \chi\biggl(\mathcal{O} \Bigl(\bigoplus_{k \geq 1} U_{2k+1}
    \otimes K\Bigr)\biggr)
    - \chi\biggl(\mathcal{O}\Bigl(\bigoplus_{k \geq 1} U_{2k}\Bigr)\biggr)\\ 
  = (g-1) \sum_{k \geq 1} \bigl(\rk(U_{2k}) + \rk(U_{2k+1})\bigr)
   + \sum_{k \geq 1} (\deg(U_{2k+1}) - \deg(U_{2k})).
\end{multline}
The stable Higgs bundle $(P_{K_{\C}},\Phi)$ represents a local minimum
of $f$ if and only if this number is zero.

Finally consider the case of reducible Higgs bundles.  It is shown in
\cite{hitchinlie} that if $\mathbb{H}^1(C^{\bullet}_{-}) = 0$ then we
continue to have a local minimum and, on the other hand, if there is
an element of $\mathbb{H}^1(C^{\bullet}_{-})$ on which $\ddot{f}$ is
negative and which is tangent to a smooth family of deformations of
the Higgs bundle, then this does not represent a local minimum.

\subsection{A theorem of Laumon}
\label{sec:laumon}
We conclude this section by a digression to the moduli space of flat
$G_{\C}$ bundles.  This space has a holomorphic symplectic form at its
smooth points (it is hyper-K{\"a}hler).  Of course the theory outlined
above applies in this case.  In particular, if the moduli space is
smooth (for example the moduli space of stable Higgs vector bundles
with rank and degree co-prime) one can consider the Morse flow on it.
At a critical point, we again have the decomposition $\Ad P_{\C} =
\bigoplus_{m}U_m$.  The subspace $T_{\leq 0}$ of the tangent space to
the moduli space on which the Hessian is less than or equal to zero is
$\mathbb{H}^1$ of the following complex of sheaves
\begin{displaymath}
  \mathcal{O}\Bigl(\bigoplus_{m \geq 0} U_{m}\Bigr)
  \xrightarrow{\ad(\Phi)} \mathcal{O}
  \Bigl(\bigoplus_{m \geq 1} U_{m} \otimes K\Bigr).
\end{displaymath}
The dimension of this is given by Riemann-Roch as
\begin{align*}
  \dim_{\C} T_{\leq 0}
  &= \chi \biggl(\mathcal{O} \Bigl(\bigoplus_{m \geq 1} U_{m} \otimes
    K\Bigr)\biggr)
    - \chi \biggl(\mathcal{O}\Bigl(\bigoplus_{m \geq 0} U_{m}\Bigr)\biggr)\\ 
  &= (g-1)\Bigl(\rk(U_0)
    + \sum_{m \geq 1} 2\rk(U_m)\Bigr) - \deg(U_0).
\end{align*}
But from \eqref{eq:UmtoU-m*} we have $\deg(U_0)=0$ and $\rk(U_0) +
\sum_{m \geq 1} 2\rk(U_m) = \dim_{\C} \lie{g}_{\C}$ so that
\begin{displaymath}
  \dim_{\C} T_{\leq 0} = \dim_{\C}(\lie{g}_{\C}) =
  \half\dim_{\C}\mathcal{M}_{G_{\C}}. 
\end{displaymath}
It was pointed out by Hausel \cite{hausel} that this fact, together with
his theorem that the downwards Morse flow coincides with the nilpotent
cone (the pre-image of $0$ under the Hitchin map),
implies a theorem of Laumon \cite[Th.\ 3.1]{laumon} in this context:
The nilpotent cone  in
$\mathcal{M}$ is a Lagrangian subvariety with respect to the
holomorphic symplectic form on $\mathcal{M}$.

\section{Milnor-Wood inequalities}
\label{sec:milnorwood}

For any $G$ there is a locally constant obstruction map
\begin{displaymath}
  o_2 \colon \Hom(\pi_1(\Sigma);G) \to H^2(\Sigma;\pi_1(G)).
\end{displaymath}
Note that in both cases, $G=\SU(n,n)$ and $G=\Sp(2n,\R)$, $\pi_1(G) =
\Z$ so that we have an integer valued function.  In the case of
representations in $\Sp(2n,\R)$, we have $o_2(\rho) = c_1(V)$, where
$V$ is the vector bundle appearing in the decomposition
\eqref{eq:sphiggsbundles} of the Higgs bundle associated to $\rho$.
In the case of $\SU(n,n)$ representations, we have $o_2(\rho) =
c_1(V)$, where $V$ is the vector bundle appearing in
\eqref{eq:suhiggsbundles}.  In both cases we thus have an integer
valued function $d = \deg(V) = \langle c_1(V), [\Sigma] \rangle$ whose
fibres are unions of connected components.  There is an outer
automorphism of $\mathcal{M}_G$ given by exchanging $V$ with $V^*$ (in
the $\Sp(2n,\R)$ case), or exchanging $V$ and $V'$ (in the $\SU(n,n)$
case).  Thus, in both cases we have an isomorphism between
$o_2^{-1}(d)$ and $o_2^{-1}(-d)$, and it therefore suffices to
consider the case $d \geq 0$, whenever convenient.

It is well known that there are bounds on the possible values of
characteristic numbers of flat bundles, known as Milnor-Wood
inequalities and, using Higgs bundles, we shall prove one for flat
$\SU(n,n)$- and $\Sp(2n,\R)$-bundles.  The original inequality proved
by J.\ Milnor \cite{milnor} concerns $\SL(2,\R)$-bundles, while J.\ 
Wood \cite{wood} considered $\SU(1,1)$-bundles.  J.\ Dupont
\cite{dupont} found a bound for any semi-simple group with finite
centre, however, the inequality of \propref{prop:milwood} below for $G
= \Sp(2n,\R)$ is sharper than his.  Using ideas of Gromov, A.\ Domic
and D.\ Toledo \cite{domic-toledo} proved a general result for
mappings of a surface into manifolds covered by bounded symmetric
domains, and their work implies \propref{prop:milwood} below.  Hitchin
obtained a proof in the case of flat reductive $\SL(2,\R)$-bundles,
using Higgs bundles, in \cite{hitchin}, and we obtain our inequality
in a similar way.  The reason why we include the proof here is, that
it gives crucial extra information about the poly-stable Higgs bundles
of the form \eqref{eq:suhiggsbundles} and \eqref{eq:sphiggsbundles}
(see \propref{prop:iso}).

\begin{prop}
\label{prop:milwood}
        Let $\rho$ be a reductive representation of $\pi_1(\Sigma)$ in
        $\Sp(2n,\R)$ or $\SU(n,n)$.  Then the characteristic number $d
        = \langle o_2(\rho),[\Sigma] \rangle$
        satisfies the inequality 
        \begin{displaymath} 
         \abs{d} \leq n(g-1).  
        \end{displaymath}
\end{prop}
\begin{proof}
We give the proof in case of $\Sp(2n,\R)$-representations, the
$\SU(n,n)$ case being completely analogous.

Let $(E,\Phi)$ be the poly-stable Higgs bundle of the form
\eqref{eq:sphiggsbundles} corresponding to $\rho$, as we already
noticed $d = \deg(V)$.  Without loss of generality we can assume that
$d>0$. In this case $c\neq 0$, as otherwise $V$ would be
$\Phi$-invariant, and therefore violate the stability condition.  Let
$U$ be the subbundle of $V^*$, such that $U \otimes K$ is the vector
bundle generated by the image of $c$.  Similarly, let $U' \subset V$
be the subbundle, which is generated by the kernel of $c$.  Then the
bundles $U'$ and $V \oplus U$ are both $\Phi$-invariant. We therefore
get the following inequalities from semi-stability of $(E,\Phi)$:
\begin{align}
        \deg(U') &\leq 0 \label{ineqone} \\ 
        d+\deg(U) &\leq 0. \label{ineqtwo} 
\end{align} 
Note that these inequalities also hold in the case when $U' = 0$ and
$U = V^*$.
Next, we note that $c$ induces a
non-trivial global section of the linebundle
\begin{displaymath} 
        \det(V/U')^{-1} \otimes \det(U \otimes K),
\end{displaymath} 
which therefore has positive degree, i.e.
\begin{equation}
\label{degline} 
        \deg(U')-d+\deg(U)+(2g-2)\rk(c) \geq 0,
\end{equation}
where $\rk(c) = \rk(U)$ is the generic rank of $c$.
Combining this with the inequalities
(\ref{ineqone}) and (\ref{ineqtwo}), we obtain
\begin{equation}
\label{eq:degreeVleq}
        d \leq (g-1) \rk(c),            
\end{equation}
so $d\leq n(g-1)$ as claimed.
\end{proof}

The above proof gives some important additional information:  from
\eqref{eq:degreeVleq} it follows that $\rk(c) = n$ for $d > (n-1)(g-1)$.
In particular, in the
extremal case $d=n(g-1)$, we have $\rk(c)=n$, and furthermore equality
holds in (\ref{degline}). Hence, $\det(c)$ is a non-zero section of a
linebundle of degree $0$, and we conclude that $c$ is an isomorphism.
We thus have the following proposition:

\begin{prop}\label{prop:iso}
  Let $(E,\Phi)$ be the poly-stable Higgs bundle of the form
  \eqref{eq:sphiggsbundles} corresponding to a reductive
  representation of $\pi_1(\Sigma)$ in $\Sp(2n,\R)$.  If
  $\deg(V)=n(g-1)$ then $c \colon V \to V^* \otimes K$ is an
  isomorphism.

  Let $(E,\Phi)$ be the poly-stable Higgs bundle of the form
  \eqref{eq:suhiggsbundles} corresponding to a reductive
  representation of $\pi_1(\Sigma)$ in $\SU(n,n)$.  If
  $\deg(V)=n(g-1)$ then $c \colon V \to V' \otimes K$ is an
  isomorphism.  
  \qed
\end{prop}

This has as a consequence that there is another discrete invariant on
$\mathcal{M}_{\Sp(2n,\R)}$ and we shall come back to this in
\secref{sec:spcomp}.

\section{Minima of $f$}
\label{sec:minima}

In this section we determine the poly-stable Higgs bundles which
represent local minima of the function $f$ on $\mathcal{M}_{G_r}$ in
the cases $G_r = \SU(2,2)$ and $G_r = \Sp(4,\R)$.  We shall determine
which stable Higgs bundles correspond to critical points of $f$ and
then identify those which are local minima using
\eqref{eq:morseindex}.  

It will be convenient to consider the decomposition \eqref{eq:decomp},
$E = \bigoplus_m F_m$ of the Higgs bundles $(E,\Phi)$ of the form
\eqref{eq:suhiggsbundles} and \eqref{eq:sphiggsbundles}, which then
gives rise to the decomposition of the adjoint bundle: note that we
have $\Ad P_{\C} = \End(V \oplus V')_0$ (the subscript $0$
indicating traceless endomorphisms) when $G_r = \SU(n,n)$, while $\Ad
P_{\C} = \End(V) \oplus S^2V \oplus S^2V^*$ for $G_r = \Sp(2n,\R)$.

We begin by finding the minima on $\mathcal{M}_{\SU(2,2)}$ which are
stable Higgs bundles, leaving the reducible ones for later.  As noted
in \secref{sec:milnorwood}, we only need to consider Higgs bundles $E=
V \oplus V'$ with $\deg(V) \geq 0$.

\begin{prop}\label{prop:suminima}
  The stable Higgs bundles of the form \eqref{eq:suhiggsbundles} with
  $\deg(V) \geq 0$, which correspond to a local minimum of $f$ on
  $\mathcal{M}_{\SU(2,2)}$ are the ones which have
  $b = 0$, $c \neq 0$, and $\deg(V) > 0$.
\end{prop}
\begin{proof}
  Let $(E,\Phi)$ be a Higgs bundle of the form
  \eqref{eq:suhiggsbundles} which represents a critical point of $f$.
  $E$ comes from the standard representation of
  $\mathrm{S}(\U(2)\times\U(2))$ on $\C^2 \oplus \C^2$ and the
  infinitesimal gauge transformation $\psi$ which produces the
  decomposition $E=\bigoplus_m F_m$ is fibrewise in
  $\lie{s}(\lie{u}(2)\times\lie{u}(2))$.  Hence each of the bundles
  $F_m$ is of the form $F_m = F_{m1} \oplus F_{m2}$ where $F_{mi} =
  V_i \cap F_m \subseteq V_i$ for $i=1,2$.  We claim that either
  $F_{m1}$ or $F_{m2}$ must be zero (unless $E=U_0$).  To see this,
  let $m_0$ be the smallest $m$ such that $F_{m1}$ and $F_{m2}$ are
  both non-zero.  Then $\Phi(U_{m_0-1})$ is contained in either $V$
  or $V'$, since the same is true for $U_{m_0-1}$ and $\Phi$
  interchanges $V$ and $V'$.  Without loss of generality we may
  suppose that $\Phi(U_{m_0-1}) \subseteq V$.  Then each of the
  bundles
  \begin{displaymath}
    \bigoplus_{m < m_0} F_m \oplus
    \Bigl(V \cap \bigoplus_{m \geq m_0} F_m\Bigr)
  \end{displaymath}
  and 
  \begin{displaymath}
    V' \cap \bigoplus_{m \geq m_0} F_m
  \end{displaymath}
  is $\Phi$-invariant, and so we have a decomposition of $(E,\Phi)$ as
  a direct sum of lower rank Higgs bundles.  This is impossible
  because $(E,\Phi)$ is stable.
  
  Let $r = (\rk(F_m))$ be the \emph{rank vector} whose entries are the
  ranks of the bundles $F_m$.  We analyze the possibilities for $r$
  case by case.

  \emph{1st case: $r=(1,1,1,1)$.}  Note that $0= \tr(\psi) =
  \ii\sum m\rk(F_m)$.  In this case we therefore have $E = F_{-3/2}
  \oplus F_{-1/2} \oplus F_{1/2} \oplus F_{3/2}$, where each $F$ is a
  linebundle.  Hence the decomposition \eqref{eq:decomp} of $\Ad P_{\C}$
  is of the form
  \begin{displaymath}
    \Ad P_{\C} = U_{-3} \oplus \cdots \oplus U_{3}, 
  \end{displaymath}
  where
  \begin{align*}
    U_2 &= \Hom(F_{-3/2},F_{1/2}) \oplus \Hom(F_{-1/2},F_{3/2}), \\
    U_3 &= \Hom(F_{-3/2},F_{3/2}).
  \end{align*}
  The formula \eqref{eq:morseindex} for the Morse index then takes the
  form
  \begin{displaymath}
  \begin{split}
    \dim_{\C} \mathbb{H}^1(C^{\bullet}_{-})
    &= 3(g-1) + \deg(F_{3/2}) - \deg(F_{-3/2}) \\
    & \quad - \bigl(\deg(F_{1/2}) - \deg(F_{-3/2}) +
    \deg(F_{3/2}) - \deg(F_{-1/2}) \bigr) \\
    &= 3(g-1) + \deg(F_{-1/2}) - \deg(F_{1/2}).
  \end{split}
  \end{displaymath}
  Now we note that $F_{1/2} \oplus F_{3/2}$ is a $\Phi$-invariant
  subbundle of $E$ and thus, by stability, $\deg(F_{1/2}) +
  \deg(F_{3/2}) < 0$.  Combining this with the  above result we get
  \begin{displaymath}
    \dim_{\C} \mathbb{H}^1(C^{\bullet}_{-})
    > 3(g-1) + \deg(F_{-1/2}) + \deg(F_{3/2}).
  \end{displaymath}
  But since $\Phi$ interchanges $V$ and $V'$ we must have $V =
  F_{-3/2} \oplus F_{1/2}$ and $V' = F_{3/2} \oplus F_{-1/2}$ or
  vice-versa.  Therefore $\abs{d} = \abs{\deg(V)} = \abs{\deg(V')}
  = \abs{\deg(F_{-1/2}) + \deg(F_{3/2})}$.  Combining this with the
  above inequality we get
  \begin{displaymath}
    \dim_{\C} \mathbb{H}^1(C^{\bullet}_{-})
    > 3(g-1) - \abs{d} \geq 0,
  \end{displaymath}
  where the last inequality comes from the Milnor-Wood inequality
  $\abs{d} \leq 2(g-1)$ of \propref{prop:milwood}.  We conclude that a
  critical point of this type always has strictly positive Morse index
  and hence it cannot be a local minimum of $f$.

  \emph{2nd case: $r=(1,2,1)$.}  Again using $\sum m \rk(F_m) = 0$ we
  see that in this case $E = F_{-1} \oplus F_{0} \oplus F_{1}$.  We
  then have $\Ad P_{\C} = U_{-2} \oplus \cdots \oplus U_2$ and
  \begin{displaymath}
    U_2 = \Hom(F_{-1},F_{1}),
  \end{displaymath}
  and so, from \eqref{eq:morseindex}, we get
  \begin{align*}
    \dim_{\C} \mathbb{H}^1(C^{\bullet}_{-})
    &= g-1 - \bigl( \deg(F_1) - \deg(F_{-1})\bigr)  \\
    &= g-1 - \bigl(2\deg(F_1) + \deg(F_{0})\bigr),
  \end{align*}
  where the second equality is due to the fact that $\deg(E) = 0$.
  Since $F_0 \oplus F_1$ and $F_1$ are $\Phi$-invariant we get from
  stability that $\deg(F_0) + \deg(F_1) < 0$ and $\deg(F_1)<0$.  Hence
  $2\deg(F_1) + \deg(F_{0}) <0 $, which shows that $\dim_{\C}
  \mathbb{H}^1(C^{\bullet}_{-}) > 0$.  Therefore a critical point of
  this type cannot be a minimum of $f$ either.

  \emph{3rd case: $r=(1,1,2)$ (or $r=(2,1,1)$).}  In this case $E =
  F_{m_1} \oplus F_{m_2} \oplus F_{m_3}$ where $V = F_{m_1} \oplus
  F_{m_2}$ and $V' = F_{m_3}$ (or vice-versa).  Since $\Phi$
  interchanges $V$ and $V'$, it follows that $\Phi_{\vert F_{m_1}}
  = 0$ and so, $(E,\Phi)$ is reducible.  Thus this case cannot occur.
  The case $r=(2,1,1)$ is analogous.
  
  \emph{4th case: $r=(2,2)$.}  In this case $E = F_{-1/2} \oplus
  F_{1/2}$.  Then $U_m = 0$ for $m \geq 2$ and hence we see from
  \eqref{eq:morseindex} that these critical points are local minima of
  $f$.  Clearly $F_{-1/2} = V$ and $F_{1/2} = V'$, or vice-versa.
  If $V = F_{1/2}$ it would be $\Phi$-invariant and so $\deg(V) <
  0$ which is absurd.  Thus, in fact, $V = F_{-1/2}$ and $V' =
  F_{1/2}$.  In the notation of \eqref{eq:suhiggsbundles} this means
  that $c = \Phi_{\vert V}$ and $b = \Phi_{\vert V'} = 0$.  This
  gives the minima with $\deg(V) > 0$.  Finally, note that if
  $\deg(V) = \deg(V') = 0$ then either $V$ or $V'$ is a
  $\Phi$-invariant subbundle which violates stability.  Thus there
  are no stable Higgs bundles with $\deg(V) = 0$ which are local
  minima of $f$.
\end{proof}

The fact that there is another discrete invariant for flat
$\Sp(2n,\R)$-bundles (cf.\ \secref{sec:spcomp}) is reflected in the
difference between the previous and the following result.

\begin{prop}\label{prop:spminima}
  The stable Higgs bundles of the form
  \eqref{eq:sphiggsbundles} with
  $\deg(V) \geq 0$, which correspond to a local minimum of $f$ on
  $\mathcal{M}_{\Sp(4,\R)}$ are the ones which have
  \begin{enumerate}
  \item $b = 0$, $c \neq 0$, and $\deg(V) > 0$.
  \item $\deg(V) = 2g-2$, $V=L_1 \oplus L_2$, and $\Phi$ of the form
    \begin{displaymath}
      \begin{pmatrix}
        0 & 0         & 0         & \tilde{c} \\
        0 & 0         & \tilde{c} & 0         \\
        0 & 0         & 0         & 0         \\
        0 & \tilde{b} & 0         & 0
      \end{pmatrix}
    \end{displaymath}
    with respect to the decomposition $E = V \oplus V^* = L_1 \oplus
    L_2 \oplus L_1^{-1} \oplus L_2^{-1}$.
  \end{enumerate}
\end{prop}
\begin{proof}
  This is analogous to the proof of
  \propref{prop:suminima}.  However, in this case the infinitesimal
  gauge transformation $\psi$ which produces the decomposition $E =
  \bigoplus F_m$ of the Higgs bundle of the form
  \eqref{eq:sphiggsbundles} belongs to $\Omega^0(\Sigma;\Ad P_K)$,
  that is, it is fibrewise in $\lie{u}(2)$.  Thus there are only two
  possibilities: either $V = F_{-1/2}$ and $V^* = F_{1/2}$ with $\Phi
  \colon V \to V^* \otimes K$, that is, $b=0$ (here we are using that
  $\deg(V) \geq 0$.  These Higgs bundles are seen to be minima as
  before.  The other possibility is that $V = F_{m_1} \oplus F_{m_2}$
  and $V^* = F_{-m_1} \oplus F_{-m_2}$, where $F_{-m} = F_{m}^*$.
  Note that either $(m_1,m_2) = (-3/2,1/2)$ or $(m_1,m_2) =
  (-1/2,3/2)$.  In
  this case the decomposition \eqref{eq:decomp} has the form 
  \begin{displaymath}
    \Ad P_{\C} = U_{-3} \oplus \cdots \oplus U_{3}, 
  \end{displaymath}
  where
  \begin{align*}
    U_2 &= \Hom(F_{-3/2},F_{1/2}) \cong \Hom(F_{-1/2},F_{3/2}), \\
    U_3 &= \Hom(F_{-3/2},F_{3/2}).
  \end{align*}
  From \eqref{eq:morseindex} we therefore get the Morse index
  \begin{displaymath}
  \begin{split}
    \dim_{\C} \mathbb{H}^1(C^{\bullet}_{-})
    &= 2(g-1) + \deg(F_{3/2}) - \deg(F_{-3/2}) \\
    & \quad   - \bigl(\deg(F_{1/2}) - \deg(F_{-3/2})\bigr) \\
    &= 2(g-1) - \bigl(\deg(F_{-3/2}) + \deg(F_{1/2})\bigr) \\
    &= 2(g-1) \pm \deg(V).
  \end{split}
  \end{displaymath}
  Thus we cannot have a minimum unless $\deg(V) = 2g-2$, and in this
  case, from \propref{prop:iso}, we have $V = F_{-3/2} \oplus F_{1/2}$
  since otherwise $c$ would not be of rank $2$.  This gives the second
  case of the proposition.
\end{proof}

It remains to identify the local minima of $f$ which are not stable
Higgs bundles.

\begin{prop}
  \label{prop:suminima2}
  The reducible Higgs bundles of the form \eqref{eq:suhiggsbundles}
  with $\deg(V) \geq 0$ which correspond to a local minimum of $f$ on
  $\mathcal{M}_{\SU(2,2)}$ either have $\Phi = 0$ and $\deg(V) =
  \deg(V') = 0$, or, if $\Phi \neq 0$, they are direct sums of rank
  $2$ Higgs bundles 
  $(E_1,\Phi_1)$ and $(E_2,\Phi_2)$, where $E_i = L_i \oplus
  L_i'$, $L_i$ a line-bundle with $\deg(L_i) \geq 0$ and $\Phi_i
  \colon L_i \to L_i' \otimes K$.  If $\Phi_i \neq 0$ then $\deg(L_i)
  > 0$.
\end{prop}

\begin{proof}
  Let $(E,\Phi)$ be a reducible Higgs bundle of the form
  \eqref{eq:suhiggsbundles} which is a local minimum of $f$.  Consider
  $\mathcal{M}_{\SU(2,2)}$ as the space of solutions $(A,\Phi)$ to
  Hitchin's equations modulo $\mathrm{S}(\U(2)\times \U(2))$ gauge
  equivalence.

  First consider the case $\Phi=0$.  Then, by poly-stability, $\deg(V)
  = \deg(V') = 0$,  and $V$ and $V'$ are poly-stable vector bundles.  On
  the other hand, it is clear that such Higgs bundles are, in fact,
  reducible (absolute) minima of $f$.  This gives the first case of
  the proposition.

  Suppose now that $\Phi \neq 0$.  The possible reductions of
  structure group are the following.
  
  \textit{Reduction to $\mathrm{S}\bigl(\bigl(\U(1) \times \U(1)\bigr)
    \times \bigl(\U(1) \times \U(1)\bigr)\bigr)$}.  In this case we
  have $V = L_1 \oplus L_2$ for linebundles $L_1$ and $L_2$, while
  $V'= L_1'\oplus L_2'$.  Thus $(E,\Phi)$ is the direct sum of two
  Higgs bundles $(E_1,\Phi_1)$ and $(E_2,\Phi_2)$, where $E_i = L_i
  \oplus L_i'$, $L_1L_1'L_2L_2'= \mathcal{O}$, and the Higgs field
  $\Phi_i$ has zeros along the diagonal.  Note also that $\deg(E_i) =
  0$ by poly-stability of $(E,\Phi)$.  Each of the bundles
  $(E_i,\Phi_i)$ is a minimum on the moduli space of rank $2$ Higgs
  bundles of this form and hence of the form \eqref{eq:decomp}, in
  other words all components of $\Phi_i$ are zero, except one
  off-diagonal entry. (cf.\ Hitchin \cite{hitchin}, Section 10).
  
  There are now two cases to consider.  The first case is when $\Phi$
  is zero on one of the bundles $V$ and $V'$; since $\deg(V) \geq 0$
  we must have $\Phi \colon V \to
  V' \otimes K$.
  In other words, $\Phi$ is of the form 
  \begin{displaymath}
    \begin{pmatrix}
      0 & 0 \\
      c & 0
    \end{pmatrix}
  \end{displaymath}
  with respect to the decomposition $E = V \oplus V'$.  Thus
  $(E,\Phi)$ is of the form considered in \propref{prop:suminima}.  As
  in the proof of that proposition one sees that there is no subspace
  of the Zariski tangent space with negative weights and, therefore,
  these Higgs bundles represent local minima of $f$.  This case
  includes the case of one of the $\Phi_i$ being equal to zero.  Note
  that, if $\deg(L_i) = 0$ then it follows from \eqref{eq:vanishing}
  and the remark following it that $\Phi_i = 0$.  This case gives the
  remaining local minima of the statement of the proposition.
  
  The other case is when $\Phi$ is non-zero on both $V$ and $V'$, say
  that $\Phi_1 \colon L_1 \to L_1'\otimes K$ and $\Phi_2 \colon L_2'
  \to L_2\otimes K$.  By stability, and since $\Phi \neq 0$, we then
  have $\deg(L_1') < 0$ and $\deg(L_2) < 0$, and so $\deg(L_1) > 0$
  and $\deg(L_2') > 0$.  We shall show that in this case $(E,\Phi)$ is
  not a local minimum of $f$.  Let the infinitesimal gauge
  transformation $\psi$ which produces the decomposition $E_i = L_i
  \oplus L_i'$ of \eqref{eq:decomp} have weights $m_i$ on $L_i$, and
  weights $m_i'$ on $L_i'$.  We then have the following equations
  relating these numbers:
  \begin{align*}
    m_1' &= m_1  + 1, \\
    m_2  &= m_2' + 1, \\
  \end{align*}
  From these equations it follows that $(m_2 - m_1) + (m_1' - m_2') =
  2$ and hence, either $m_2 - m_1 \geq 1$ or $m_1'-m_2'\geq 1$.  For
  definiteness suppose that $m_2 - m_1 \geq 1$ (the other case is
  entirely similar).  This means that
  \begin{displaymath}
    \Hom(L_1,L_2) \subseteq \Ad_{\lie{k}_{\C}} P_K = \bigl(\End(V)
    \oplus \End(V') \bigr)_0
  \end{displaymath}
  has weight $\geq 1$ and that this is a subspace of the highest
  weight space of $\psi$.  Note that $\ad(\Phi)$ is zero restricted to
  the highest weight space and so $H^1(\Sigma;\Hom(L_1,L_2))$ gives a
  subspace of $\mathbb{H}^1(C^{\bullet}_{-})$ on which $\ddot{f}$ is
  negative.  But since $\deg(L_1) \geq 0$ and $\deg(L_2) < 0$ we have
  $H^0(\Sigma;\Hom(L_1,L_2)) = 0$ and therefore, from Riemann-Roch,
  $H^1(\Sigma;\Hom(L_1,L_2)) \neq 0$.  It only remains to find a
  smooth family of Higgs bundles in $\mathcal{M}_{\SU(2,2)}$ to which
  an element in $H^1(\Sigma;\Hom(L_1,L_2))$ is tangent (cf.\ Section
  \ref{sec:morseindices}).  By hypothesis $(E,\Phi)$ is the direct sum
  of the stable Higgs bundles $(E_1,\Phi_1)$ and $(E_2,\Phi_2)$.  All
  extensions
  \begin{displaymath}
    0 \to E_2 \to E \to E_1 \to 0
  \end{displaymath}
  are parametrized by $H^1(\Sigma; \Hom(E_1,E_2))$ so, in particular,
  $t \in H^1(\Sigma; \Hom(L_1,L_2))$ defines an extension
  \begin{displaymath}
    0 \to E_2 \to E_t \to E_1 \to 0
  \end{displaymath}
  which is non-trivial if $t \neq 0$.  Note that $E = V_t \oplus V'$,
  where $V_t$ is the non-trivial extension 
  \begin{displaymath}
    0 \to L_2 \to V_t \to L_1 \to 0
  \end{displaymath}
  defined by $t$.  We define a Higgs field 
  \begin{math}
    \Phi =
    \left(
    \begin{smallmatrix}
      0 & b_t \\
      c_t & 0
    \end{smallmatrix}
    \right)
  \end{math}
  on $E_t$ of the appropriate form in the following way.  To define
  $b_t \colon V'\to V_t \otimes K$ we use the composition
  \begin{displaymath}
    V' \xrightarrow{b} L_2 \otimes K \to V_t \otimes K,
  \end{displaymath}
  while to define
  $c_t \colon V_t\to V' \otimes K$ we use the composition
  \begin{displaymath}
    V_t \to L_1 \xrightarrow{c} V' \otimes K.
  \end{displaymath}
  For $t \neq 0$ the Higgs bundle $(E_t,\Phi_t)$ is a non-trivial
  extension of stable Higgs bundles and therefore stable.  For $\alpha
  \in \C$, the family $(E_{\alpha t},\Phi_{\alpha t})$ is thus a
  smooth family of Higgs bundles in $\mathcal{M}_{\SU(2,2)}$ to which
  $t \in H^1(\Sigma;\Hom(L_1,L_2))$ is tangent.
  
  \textit{Reduction to $\mathrm{S}\bigl(\bigl(\U(1) \times \U(1)\bigr)
  \times \U(2)\bigr)$}.
  In this case we have a decomposition of $(E,\Phi)$ as a direct sum
  of Higgs bundles $(E_1,\Phi_1)$ and $(E_2,\Phi_2$), where $E_1 = V
  \oplus L_1$ and $E_2 = L_2$ with $L_1$ and $L_2$ linebundles.
  Again $(E_1,\Phi_1)$ and $(E_2,\Phi_2$) represent local minima on
  lower rank moduli spaces and so, $\Phi_2 = 0$.  If $(E_1,\Phi_1)$ is
  reducible we are back in one of the previous cases so we may assume
  that $(E_1,\Phi_1)$ is stable.  Since we are at a minimum it must be
  of the form \eqref{eq:decomp} and again there are several
  possibilities.  If $\Phi_1$ is zero on either $V$ or $L_1$ then it
  is zero on either $V$ or $V'$ and $(E,\Phi)$ is a minimum as above.
  Thus the only case that remains is when $V = F_{-1} \oplus F_{1}$
  and $\Phi \colon F_{-1} \to L_1 \otimes K$, and $\Phi \colon L_1 \to
  F_1 \otimes K$.  The weights of the infinitesimal gauge
  transformation producing this decomposition are $-1$, $0$, and $1$
  on $F_{-1}$, $L_1$, and $F_1$, respectively.  As above one sees that
  $H^1(\Sigma ; \Hom(F_{-1}, F_1))$ gives a subspace of
  $\mathbb{H}^1(C^{\bullet}_{-})$ on which $\ddot{f}$ is negative and
  that $t \in H^1(\Sigma ; \Hom(F_{-1}, F_1))$ is tangent to a smooth
  family of Higgs bundles, so that $(E,\Phi)$ does not represent a
  local minimum.  We omit the details.
  
  \textit{Reduction to $\mathrm{S}\bigl(\U(2) \times \bigl(\U(1) \times
  \U(1)\bigr)\bigr)$}.  This case is analogous to the previous one. 
\end{proof}

In an analogous manner one can prove the following proposition.

\begin{prop}
  \label{prop:spminima2}
  The reducible Higgs bundles of the form \eqref{eq:sphiggsbundles}
  with $\deg(V) \geq 0$ which correspond to a local minimum of $f$ on
  $\mathcal{M}_{\Sp(2,\R)}$ either have $\Phi = 0$ and $\deg(V) = 0$,
  or, if $\Phi \neq 0$, they are direct sums of rank $2$ Higgs bundles
  $(E_1,\Phi_1)$ and $(E_2,\Phi_2)$, where $E_i = L_i \oplus
  L_i^{-1}$, $L_i$ a line-bundle with $\deg(L_i) \geq 0$ and $\Phi_i
  \colon L_i \to L_i^{-1} \otimes K$.  If $\Phi_i \neq 0$ then
  $\deg(L_i) > 0$.\qed
\end{prop}

\section{Connected Components}
\label{sec:count}

\subsection{Components of $\mathcal{M}_{\SU(2,2)}$}
\label{sec:sucomp}

In this section we consider the connected
components of $\mathcal{M}_{\SU(2,2)}$.  Using \propref{prop:milwood}
we can write 
\begin{displaymath}
  \mathcal{M}_{\SU(2,2)} = \mathcal{M}_{-(2g-2)} \cup \cdots \cup
  \mathcal{M}_{2g-2},
\end{displaymath}
where $\mathcal{M}_d$, the subspace of Higgs bundles of the form
\ref{eq:suhiggsbundles} with $\deg(V) = d$, is a union of connected
components.  Denote the subspace of local minima of $f$ on
$\mathcal{M}_d$ by $\mathcal{N}_d$.  Note that, since $f$ is proper,
connectedness of $\mathcal{N}_d$ implies connectedness of
$\mathcal{M}_d$.  As noted in \secref{sec:milnorwood}, we can without
loss of generality assume that $d \geq 0$.  The results of the previous
section then give the following identification of $\mathcal{N}_d$.

\begin{prop}
  \label{prop:allsuminima}
  The subspace of local minima of $f$ on $\mathcal{M}_d$,
  $\mathcal{N}_d$, is the space of poly-stable Higgs bundles of the
  form \ref{eq:suhiggsbundles} with $b = 0$.
\end{prop}

\begin{proof}
  Immediate from Propositions \ref{prop:suminima} and
  \ref{prop:suminima2}. 
\end{proof}

We can use this to identify $\mathcal{N}_d$ with a moduli space of
triples, as studied by Bradlow and Garc{\'\i}a-Prada
\cite{bradlow-garcia-prada, garcia-prada}.  Denote the moduli space of
stable triples $(V,\tilde{V'},\phi)$ (where $\phi \in
H^0(\Sigma;\Hom(V,\tilde{V'}))$, $\deg(V) = d$, $\deg(\tilde{V'}) =
2g-2-d$, and $\rk(V) = \rk(\tilde{V'}) = 2$) by
$\mathcal{M}_{d}^{\mathrm{triples}}$ (cf.\ 
\secref{sec:q-bundles-triples}).  To each such triple we associate a
Higgs bundle $ (V \oplus V', \left(
\begin{smallmatrix}
  0 & 0 \\
  \phi & 0
\end{smallmatrix}
\right))
$,
where $V' = \tilde{V'} \otimes K^{-1}$.
The following theorem is then an immediate consequence of 
\thmref{thm:quiver}.
\begin{thm}
$\mathcal{N}_d$ is isomorphic to the
fibre, over the trivial bundle $\mathcal{O}$, of the map
\begin{align*}
  \mathcal{M}_{d}^{\mathrm{triples}} &\to \Jac(\Sigma) \\
  (V,\tilde{V'},\phi) &\mapsto \Lambda^2(V)\otimes
  \Lambda^2(\tilde{V'}) \otimes K^{-2}.
\end{align*}
\qed
\end{thm}
Thus information about connectedness of moduli spaces of stable
triples would give information about connectedness of the
$\mathcal{M}_d$ and we hope to come back to this on a later occasion.
At present we can prove the following theorem.

\begin{thm}
  The subspaces $\mathcal{M}_d$ of $\mathcal{M}_{\SU(2,2)}$ are
  connected for $d=0$ and $d=\pm (2g-2)$.
\end{thm}

\begin{proof}
  First consider the case $d=0$.  To see that $\mathcal{N}_0$ is
  connected, consider the continuous map
  \begin{align*}
    N_0 \times N_0 \times \Jac(\Sigma) &\to \mathcal{N}_0 \\
    (E,E',L) &\mapsto ((E \otimes L) \oplus (E' \otimes L^{-1}), 0),
  \end{align*}
  where $N_0$ denotes the moduli space of rank $2$
  poly-stable vector bundles with fixed trivial determinant bundle.
  From \propref{prop:allsuminima} we
  see that this is surjective and, since $N_0$ and $\Jac(\Sigma)$ are
  connected, that $\mathcal{N}_0$ is connected.

  Next consider the case $d = 2g-2$ (as already noticed, this also
  takes care of the case $d= -(2g-2)$).  From Propositions
  \ref{prop:iso} and \ref{prop:allsuminima} we see that
  $\mathcal{N}_{2g-2}$ is isomorphic to the moduli space of
  rank $2$, degree $2g-2$ vector bundles with fixed determinant, which
  is known to be connected.
\end{proof}

\subsection{Components of $\mathcal{M}_{\Sp(4,\R)}$}
\label{sec:spcomp}

In this section we consider the connected
components of $\mathcal{M}_{\Sp(4,\R)}$.  Again using
\propref{prop:milwood} 
we can write 
\begin{displaymath}
  \mathcal{M}_{\Sp(4,\R)} = \mathcal{M}_{-(2g-2)} \cup \cdots \cup
  \mathcal{M}_{2g-2},
\end{displaymath}
where $\mathcal{M}_d$, the subspace of Higgs bundles of the form
\ref{eq:sphiggsbundles} with $\deg(V) = d$, is a union of connected
components.  Again we denote the subspace of local minima of $f$ on
$\mathcal{M}_d$ by $\mathcal{N}_d$, and connectedness of
$\mathcal{N}_d$ implies connectedness of $\mathcal{M}_d$.  

We can also identify $\mathcal{N}_d$ with a
moduli space of triples, using \thmref{thm:quiver}, as follows.

\begin{thm}
For $-(2g-2) \leq d \leq 2g-2$, $\mathcal{N}_d$ is the isomorphic to the
fixed point set of the involution on the moduli space
$\mathcal{M}_{d}^{\mathrm{triples}}$ of poly-stable triples (as defined
in the previous section), defined by
\begin{align*}
  \mathcal{M}_{d}^{\mathrm{triples}} &\to \mathcal{M}_{d}^{\mathrm{triples}} \\
  (V,\tilde{V'},\phi) &\mapsto \Lambda^2(V)\otimes
  \Lambda^2(\tilde{V'}) \otimes K^{-2}.
\end{align*}
\qed
\end{thm}

With regard to connectedness, we consider
the cases $\abs{d} < 2g-2$ and  $\abs{d} = 2g-2$ separately.  

\paragraph{The case $\abs{d} < 2g-2$.}
In this case everything is completely analogous to the case of
$\SU(2,2)$-bundles.  To begin with, we have the following result.

\begin{prop}
  \label{prop:allspminima}
  For $\abs{d} < 2g-2$, the subspace of local minima of $f$ on
  $\mathcal{M}_d$, $\mathcal{N}_d$, is the space of poly-stable Higgs
  bundles of the form \ref{eq:sphiggsbundles} with $b = 0$.
\end{prop}

\begin{proof}
  Follows from Propositions \ref{prop:spminima} and
  \ref{prop:spminima2}. 
\end{proof}

We have the following result about connectedness of $\mathcal{M}_0$.

\begin{thm}
  The subspace $\mathcal{M}_0$ of $\mathcal{M}_{\Sp(4,\R)}$ is
  connected.
\end{thm}

\begin{proof}
  From \propref{prop:allspminima} it follows in particular that
  $\mathcal{N}_0$ is isomorphic to the moduli space of rank
  $2$, degree $0$
  poly-stable vector bundles.  Since this space is connected, the
  result is proved.
\end{proof}

\paragraph{The case$\abs{d} = 2g-2$.}
In this case the results are
entirely different from those of $\SU(2,2)$-bundles, due to
\propref{prop:spminima}.  

Let $(E,\Phi)$ be a Higgs bundle of the form \eqref{eq:sphiggsbundles}
with $d = n(g-1)$.  
Choosing a square root $L_0$ of the canonical bundle on $\Sigma$, we 
can define a rank $n$ vector bundle $W$ by
\begin{displaymath}
        W = V \otimes L_0^{-1},
\end{displaymath}
and we can define $C \in H^0(\Sigma; S^2 W^*)$ and 
$\phi \in H^0(\Sigma; \End(W) \otimes K^2)$ by
\begin{displaymath}
        C = c \otimes 1_{L_0^{-1}},
\end{displaymath}
and
\begin{displaymath}
        \phi = (b \otimes 1_{L_0}) \circ (c \otimes 1_{L_0^{-1}}).
\end{displaymath}
Note that $\phi$ is symmetric with respect to the quadratic form $C$.

From \propref{prop:iso} we know that $c$ is an isomorphism when
$(E,\Phi)$ is poly-stable, and thus we can recover $(E,\Phi)$ from
this data.  Therefore the set of isomorphism classes of Higgs bundles
of the form \eqref{eq:sphiggsbundles} is equal to the set of
isomorphism classes of Higgs bundles
\begin{equation}
\label{eq:symhiggs}
        (W,C,\phi),
\end{equation}
where $W$ has a non-degenerate quadratic form $C$, and the Higgs field 
$\Phi$ is twisted by $K^2$ and symmetric with respect to $C$.  There
is an obvious stability condition for $(W,C,\phi)$, namely that 
\begin{equation}
\label{eq:Wstability}
        \mu(U) < \mu(W)
\end{equation}
for all $\phi$-invariant subbundles $U$ of $W$.  Next, we shall
prove that $(W,C,\phi)$ is stable, if and only if $(E,\Phi)$ is.

\begin{thm}
\label{thm:extremalcase}
The subspace $\mathcal{M}_{n(g-1)} \subset 
\mathcal{M}_{\Sp(2n,\R)}$ of Higgs bundles of the form 
\eqref{eq:sphiggsbundles}, with $d = n(g-1)$ is isomorphic to the moduli 
space of poly-stable Higgs bundles of the form \eqref{eq:symhiggs}.
\end{thm}
\begin{proof}
  We have to prove that $(E,\Phi)$ is stable if and only if
  $(W,C,\phi)$ is.  From \thmref{thm:quiver} we know that stability of
  $(E,\Phi)$ is equivalent to stability of the $Q$-bundle $\mathbf{E}
  = (\underline{E},\underline{\Phi})$. Thus, all we need to prove is
  that $\mathbf{E}$ is stable if and only if $(W,C,\phi)$ is.  Because
  stability is unaffected by tensoring with a line bundle, we can
  equally well prove that $(V,b \circ c)$ is stable.  Note, that
  $\mu(V) = g-1$.
        
  Assume $\mathbf{E}$ is a stable $Q$-bundle. Let $U \subset V$ be a
  $\phi$-invariant subbundle.  Let $U' \subset V^*$ be the subbundle
  such that $U' \otimes K$ is generically the image of $U$ under $c$.
  Then $b$ maps $U'$ to $U$, because of the $\phi$-invariance of $U$.
  Hence, $\mathbf{F} = (\{U,U'\},\{b,c\})$ defines a $Q$-subbundle of
  $\mathbf{E}$, and it follows that
  \begin{equation}
  \label{eq:stable}
     \mu(\mathbf{F}) < \mu(\mathbf{E}).
  \end{equation}
  But, as $c$ is an isomorphism
  \begin{align*}
    \mu(\mathbf{E})   &= \mu(E) \\
    &= \mu(V \oplus V \otimes K^{-1}) \\
    &= \mu(V) - (g-1),
  \end{align*}
  and similarly $\mu(\mathbf{F}) = \mu(U) - (g-1)$. Therefore $\mu(U)
  < \mu(V)$ and so, $(W,C,\phi)$ is stable.
  
  Conversely, assume that $(W,C,\phi)$ is stable. Let $\mathbf{F} =
  (\{U,U'\},\{b,c\})$ be a $Q$-subbundle of $\mathbf{E}$.  Let
  $\tilde{U} \subset V^*$ be the subbundle which is generically the
  image of $U' \otimes K$ under $c^{-1}$.  Both $U$ and $\tilde{U}$
  are $\phi$-invariant subbundles of $V$, because $\mathbf{F}$ is a
  $Q$-subbundle.  Hence, $\mu(U) < \mu(V)$ and $\mu(U') < \mu(V)$, by
  stability of $(W,C,\phi)$.  Recalling that $\mu(V) = g-1$ and
  $\mu(\tilde{U}) = \mu(U') - (2g-2)$, we get
  \begin{align}
    \mu(U)  &< g-1, \label{eq:muU} \\
    \intertext{and} \mu(U') &< -(g-1). \label{eq:muUprime}
  \end{align}
  Note also that
  \begin{equation}\label{eq:rkUprime}
    \rk(U') \geq \rk(U),
  \end{equation}
  because $c$ is an isomorphism, and the image of $U$ under $c$ is
  contained in $U' \otimes K$, by the assumption that $\mathbf{F}$ is
  a $Q$-subbundle.  Combining \eqref{eq:muU}, \eqref{eq:muUprime}, and
  \eqref{eq:rkUprime}, we get:
  \begin{align*}
    \mu(\mathbf{F}) &= \mu(U \oplus U') \\
    &= \frac{\rk(U)}{\rk(U \oplus U')} \mu(U)
    + \frac{\rk(U')}{\rk(U \oplus U')} \mu(U') \\
    &< \frac{\rk(U) - \rk(U')}{\rk(U \oplus U')} (g-1) \\
    &\leq 0.
  \end{align*}
  Of course, $\mu(\mathbf{E}) = 0$ and hence the proof is finished.
\end{proof}

The existence of the quadratic form $C$ on $W$ means that the
structure group is $\Or(n,\C)$. The maximal compact subgroup of
$\Or(n,\C)$ is $\Or(n)$ and, therefore, we have the Stiefel-Whitney
classes $w_1$ and $w_2$ as topological invariants.  We now specialize
to the case $n=2$.  The first Stiefel-Whitney class can then be seen
in holomorphic terms as follows: the quadratic form $C$ gives an
isomorphism $(\Lambda^2 W)^2 \cong \calo$; hence, $\Lambda^2 W$ gives
an element of $H^1(\Sigma; \Z / 2)$, and it is easy to see that this
element is $w_1(W)$.  It follows that $\Lambda^2 W = \calo$ if and
only if $w_1(W) = 0$.  This, in turn, is equivalent to the existence
of a reduction of structure group to $\SO(2,\C) \subset \Or(2,\C)$.
Using the identification $\C^\times \cong \SO(2,\C)$ via
\begin{displaymath}
  \lambda \mapsto
  \begin{pmatrix}
    \lambda & 0 \\
    0 & \lambda^{-1}
  \end{pmatrix},
\end{displaymath}
we see that this happens exactly when $W$ decomposes as a direct sum
\begin{displaymath}
  W = L \oplus L^{-1},
\end{displaymath}
and $C$ is of the form
\begin{displaymath}
  \begin{pmatrix}
    0 & 1 \\
    1 & 0
  \end{pmatrix}
\end{displaymath}
with respect to this decomposition.  Now it is clear that, in
this case, $w_2(W)$ is given by
\begin{displaymath}
        w_2 = c_1(L) \mod{2}.
\end{displaymath}
By interchanging $L$ with its dual if necessary, we may assume that
$\deg(L) \geq 0$.  Furthermore, when $\deg(L) > 0$, the Higgs field
$\phi$ must induce a non-zero holomorphic map
\begin{displaymath}
  L \to L^{-1} K^2,
\end{displaymath}
because otherwise $L \subset W$ would violate stability.  Hence, we have
\begin{displaymath}
  \deg(L) \leq 2g-2.
\end{displaymath}

We, therefore, have a decomposition of $\mathcal{M}_{2g-2}$ into
subspaces, each of which is a union of connected components, as
follows: 
\begin{displaymath}
  \mathcal{M}_{2g-2} = \Bigl(
                       \bigcup_{u,v}
                       \mathbf{M}^v_u
                       \Bigr)\cup
                       \Bigl(
                       \bigcup_{l=0}^{2g-2} \mathbf{M}_0^l
                       \Bigr),
\end{displaymath}
where $\mathbf{M}^v_u$ is the moduli space of poly-stable Higgs bundles
$(W,C,\phi)$ with $w_1(W) = u \in H^1(\Sigma; \Z/2)-\{0\}$
and $ w_2(W) = v \in H^2(\Sigma; \Z/2)$, and where $\mathbf{M}_0^l$
is the moduli space of poly-stable Higgs bundles
$(W,C,\phi)$ with $w_1(W) = 0$ and $\deg(L) = l$.
We shall prove that each of these subspaces is connected, except
$\mathbf{M}_0^{2g-2}$:  in this case the Higgs field $\phi$ induces a
non-zero section of the degree $0$ linebundle $L^{-2}K^2$ and thus
$L^2=K^2$.  We therefore have a further decomposition of of
$\mathbf{M}_0^{2g-2}$ into subspaces $\mathbf{M}_{0,L}^{2g-2}$,
indexed by the $2^{2g}$ square roots $L \in \Jac^{2g-2}(\Sigma)$ of
$K^2$  (note the analogy with the breaking up of $\mathcal{M}_{2g-2}$
into several connected components).

We can now state our main result, to be proved in the remaining
part of this section.

\begin{thm}
\label{thm:components}

\begin{itemize}
\item[$i)$] The spaces $\mathbf{M}_{0,L}^{2g-2}$ are connected.
\item[$ii)$] The spaces $\mathbf{M}^l_0$ are connected for $0 \leq l <
  2g-2$.
\item[$iii)$] The spaces $\mathbf{M}^v_u$ are connected.
\end{itemize}
\end{thm}

\begin{rem}
  Hitchin showed in \cite{hitchinlie} that for any split real form
  $G_r$ of a complex simple Lie group the moduli space of reductive
  representations of $\pi_1(\Sigma)$ in $G_r$ contains a connected
  component which is homeomorphic to an Euclidean space of dimension
  $(2g-2)\dim G_r$.  This component is called the Teichm{\"u}ller
  component.  The group $\Sp(4,\R)$ is a split real form of
  $\Sp(4,\C)$ so there is a Teichm{\"u}ller component in this case.
  As a matter of fact, each of the subspaces $\mathbf{M}_{0,L}^{2g-2}$
  is isomorphic to a vector space: note that $W = L \oplus L^{-1}$ is
  completely determined by $L$ and that any $(W,C,\phi)$ is stable.
  Hence $\mathbf{M}_{0,L}^{2g-2}$ is isomorphic to the space of Higgs
  fields $\phi \in H^0(\Sigma; \End(W) \otimes K^2)$ which are
  symmetric with respect to $C = \left(
    \begin{smallmatrix}
      0 & 1 \\
      1 & 0
    \end{smallmatrix}
  \right)$, that is, of the form $\phi = \left(
    \begin{smallmatrix}
      \phi_{11} & \phi_{12} \\
      \phi_{12} & \phi_{22} \\
    \end{smallmatrix}
  \right)$.  It follows that  $\mathbf{M}_{0,L}^{2g-2}$ is isomorphic
  to the vector space
  \begin{displaymath}
    H^0(\Sigma;K^2) \oplus
    H^0(\Sigma;K^2) \oplus H^0(\Sigma;K^4).
  \end{displaymath}
  Note that this proves $i)$ of the theorem.
\end{rem}

\begin{rem}
  One can see (see \cite{thesis} for details), that the subspaces
  $\mathcal{M}_d$, $\mathbf{M}^v_u$, $\mathbf{M}^l_0$, and
  $\mathbf{M}^{2g-2}_{0,L}$ are non-empty.  Therefore,
  \thmref{thm:components} shows that $\mathcal{M}_{\Sp(4,\R)}$ has at
  least $3\cdot 2^{2g} + 8g-13$ connected components.
\end{rem}

\subparagraph{Proof that the subspaces $\mathbf{M}^l_0 \subset
  \mathcal{M}_{2g-2}$ are connected.} 
Recall that any $(W,C,\phi)$ in $\mathbf{M}^l_0$ is of the form 
\begin{displaymath}
  W = L \oplus L^{-1},
\end{displaymath}
with $l = \deg(L)$ and $C$ of the form 
$
\left(
\begin{smallmatrix}
  0 & 1 \\
  1 & 0
\end{smallmatrix}
\right) $.  First, we consider the case of $l > 0$.  In this case, the
Higgs field $\phi$ must be non-zero, as otherwise the subbundle $L
\subset W$ would violate stability.  But any critical point of the
type described in $i)$ of \propref{prop:spminima} has $\phi = 0$ so,
it follows that all the critical points in $\mathbf{M}^l_0$ for $l >
0$ are of the type described in $ii)$ of \propref{prop:spminima}.  We
therefore see that the critical points correspond to Higgs bundles
$(W,C,\phi)$, which are of the form described above and where,
furthermore, $\phi$ is of the form
\begin{displaymath}
  \phi =
  \begin{pmatrix}
    0 & 0 \\
    \tilde{\phi} & 0
  \end{pmatrix},
\end{displaymath}
with $\tilde{\phi} \in H^0(\Sigma; L^{-2}K^2)$.  Using this, it is now
easy to give an explicit description of the subspace of local minima
of $f$ on $\mathbf{M}^l_0$.

\begin{prop}
The subspace of local minima $N^l_0 \subset \mathbf{M}^l_0$ fits into a
pull-back diagram
\begin{displaymath}
\begin{CD}
  N^l_0             @>>>                \Jac^l(\Sigma) \\
  @VV{\pi}V                             @VV{L \mapsto L^{-2}K^2}V \\
  S^{4g-4-2l}\Sigma @>{D \mapsto [D]}>> \Jac^{4g-4-2l}(\Sigma),
\end{CD}
\end{displaymath}
where $\pi(W,C,\phi) = (\phi)$.
\end{prop}

\begin{proof}
The only thing there is to remark is that any $(W,C,\phi)$, of the
form given above, is stable.  But, $L^{-1} \subset W$ is the
only $\phi$-invariant subbundle so, this is obvious.
\end{proof}

From this proposition, it is clear that $N^l_0$ is connected so, from
the properness of $f$, it follows that $\mathbf{M}^l_0$ is connected
for $l > 0$. 

In the case $l = 0$,  we have the following result.

\begin{prop}
  Any local minimum of $f$ on $\mathbf{M}^0_0$ has $\phi = 0$ and is,
  therefore, of the type described in $i)$ of
  \propref{prop:spminima}. 
\end{prop}

\begin{proof}
Suppose we have a critical point of the type described
in $ii)$ of \propref{prop:spminima}, with
$\phi \neq 0$.  Then, $L^{-1} \subset W$ is $\phi$-invariant and
therefore, $(W,C,\phi)$ is semi-stable, but not stable.  Since we
are considering the moduli space of poly-stable Higgs bundles, 
$(W,C,\phi)$ decomposes as a direct sum of rank $1$ Higgs bundles of
degree $0$.  The only subbundles of $W$ of rank $1$ and degree $0$
are $L$ and $L^{-1}$, and $L$ is not $\phi$-invariant so, we
conclude that this situation cannot occur.
\end{proof}

Consequently, we have the following description of the subspace of
local minima of $f$ on $\mathbf{M}^0_0$.

\begin{prop}
The subspace $N^0_0 \subset \mathbf{M}^0_0$ of local minima of $f$ is
isomorphic to the moduli space of poly-stable $(W,C)$, where $W$ is of
the form
\begin{displaymath}
  W = L \oplus L^{-1},
\end{displaymath}
for a linebundle $L$ of degree $0$, and $C$ is of the form 
$
\left(
\begin{smallmatrix}
  0 & 1 \\
  1 & 0
\end{smallmatrix}
\right)
$, with respect to this decomposition.
\end{prop}

Note that the pair $(W,C)$ decomposes into a direct sum of Higgs
linebundles exactly when $L^2 = \calo$, and it is then poly-stable,
but not stable.  All other $(W,C)$ are stable.  It follows that there
is a surjective continuous map
\begin{displaymath}
  \Jac^0(\Sigma) \to N^0_0,
\end{displaymath}
given by taking $L$ to $(W,C)$ of the form given above.  Therefore,
$N^0_0$ is connected, finishing the proof that the subspaces
$\mathbf{M}^l_0$ are connected.

\subparagraph{Proof that the subspaces $\mathbf{M}^v_u \subset
  \mathcal{M}_{2g-2}$ are connected.}  We begin by noting that the
$(W,C,\phi)$ corresponding to a critical point of the type described
in $ii)$ of \propref{prop:spminima} has $w_1(W) = 0$, thus we see that
the subspaces of local minima $N^v_u \subset \mathbf{M}^v_u$ consist
of critical points of the type described in $i)$ of
\propref{prop:spminima}.  Recall that for these $b = 0$; in terms of
the Higgs bundle $(W,C,\phi)$, this means that $\phi = 0$.  Thus,
$N^v_u$ is the moduli space of stable pairs $(W,C)$ with the given
characteristic classes.  From $\Lambda^2 W \neq \calo$, one sees
easily that any such pair is stable.

There is a connected double cover 
$\tilde{\Sigma} \overset{\pi}{\to} \Sigma$ given by 
\begin{displaymath}
w_1(W) \in H^1(\Sigma; \Z/2) = \Hom(\pi_1 \Sigma,\Z/2).
\end{displaymath}
Clearly, the pull-back of $W$ to $\tilde{\Sigma}$ is of the form
$\pi^*W = M \oplus M^{-1}$ with $ \pi^*C = 
\left(
\begin{smallmatrix}
  0 & 1 \\
  1 & 0
\end{smallmatrix}
\right)
$
and $M$ a linebundle.  Let $\tau \colon \tilde{\Sigma} \to
\tilde{\Sigma}$ be the involution interchanging the sheets of the
covering, then, clearly,
\begin{displaymath}
  \tau^*M = M^{-1}.
\end{displaymath}
Conversely, if $M$ is a linebundle on $\tilde{\Sigma}$ which satisfies
this condition, then $W = \pi_*M$ is a rank $2$ vector bundle with a
non-degenerate quadratic form $C$.  In fact $\tilde{\Sigma}$ is the
spectral curve associated to $(W,C)$ (see Hitchin \cite{hitchinduke}
and Beauville, Narasimhan and Ramanan \cite{beauville}).
Hence $N^0_u \cup N^1_u$ can be identified with the kernel of the map
\begin{displaymath}
  1 + \tau^* \colon \Jac(\tilde{\Sigma}) \to \Jac(\tilde{\Sigma}),
\end{displaymath}
where $\tilde{\Sigma}$ is the unramified double cover of $\Sigma$
given by $u \in H^1(\Sigma; \Z/2)$.

It remains to distinguish between $w_2$ being equal to 0 or 1.  When
the cover is unramified, 
the kernel of $1 + \tau^*$ splits into two components,
\begin{displaymath}
  \ker(1 + \tau^*) = P^+ \cup P^-,
\end{displaymath}
each of them a translate of the Prym variety of the covering.
It is
a classical theorem of Wirtinger, that the function $\delta \colon
P^+ \cup P^- \to \Z/2$, defined by
\begin{align*}
  \delta(M) 
    &= \dim_{\C} H^0(\tilde{\Sigma}; M \otimes \pi^* L_0) \mod 2 \\
    &= \dim_{\C} H^0(\tilde{\Sigma}; \pi_*M \otimes L_0) \mod 2,
\end{align*}
is constant on each of $P^+$ and $P^-$ and takes different values on
them.  For proofs of these facts, see Mumford \cite{mumford:theta} or
\cite{mumford:prym}. 

Now, let $F \to \Sigma$ be a real vector bundle.  Choosing a metric on
$F$, the complexification $F^c = F \otimes_{\R} \C$ acquires a
holomorphic structure and therefore, there is a $\dbar$-operator
\begin{displaymath}
  \dbar_{L_0}(F) \colon \Omega^0(\Sigma; L_0 \otimes F^c) 
                 \to \Omega^{0,1}(\Sigma; L_0 \otimes F^c).
\end{displaymath}
Atiyah \cite{atiyah:spin} shows that the function
\begin{displaymath}
  \delta_{L_0}(F) = \dim_{\C} \ker (\dbar_{L_0}(F)) \mod 2
\end{displaymath}
is independent of the choice of the metric, and that it extends to give
a group homomorphism
\begin{displaymath}
  \delta_{L_0} \colon KO(\Sigma) \to \Z/2.
\end{displaymath}
Define $\gamma \in \tKO (\Sigma)$ to be the
pull-back of the generator of $\tKO (S^2)$ under a map $\tilde{\Sigma}
\to S^2$ of degree $1$.  Atiyah \cite[Lemma (2.3)]{atiyah:spin}
shows that 
\begin{displaymath}
  \delta_{L_0}(\gamma) = 1.
\end{displaymath}
Furthermore, the total Stiefel-Whitney class gives an isomorphism
\begin{displaymath}
  w \colon \tKO(\Sigma) 
    \to \{1\} \oplus H^1(\Sigma; \Z/2) \oplus H^2(\Sigma; \Z/2)
\end{displaymath}
of the additive group $\tKO(\Sigma)$ onto the multiplicative group of
the cohomology ring $H^*(\Sigma; \Z/2)$ (see 
\cite[Remark, p.\ 54]{atiyah:spin}).  Clearly,
\begin{displaymath}
  w(\gamma) = (1,0,1),
\end{displaymath}
where we identify $H^2(\Sigma; \Z/2) = \Z/2$.  We may, therefore, think
of $\delta_{L_0}$ as a homomorphism of the multiplicative group of
$H^*(\Sigma; \Z/2)$ to $\Z/2$, which takes the value $1$ on the
element $(1,0,1)$.  Let $u \in H^1(\Sigma; \Z/2)$; then,
\begin{displaymath}
  (1,u,0) = (1,u,1) \cdot (1,0,1)
\end{displaymath}
in $H^*(\Sigma; \Z/2)$.  Therefore,
\begin{equation}
\label{eq:w2}
  \delta_{L_0}(1,u,0) = \delta_{L_0}(1,u,1) + 1.
\end{equation}

Returning to $(W,C)$ with $W = \pi_* M$ for $M \in \ker(1 + \tau^*)$,
we see that 
\begin{displaymath}
  \delta(M) = \delta_{L_0}(W^r),
\end{displaymath}
where $W^r$ is a real rank two bundle, whose complexification is $W$.
It follows from \eqref{eq:w2}, that $\delta$ takes different values
for different values of $w_2(W)$ and hence, that $w_2(W)$ determines
whether $M$ lies in $P^+$ or $P^-$.

From this discussion, we obtain the following explicit description of
the subvariety $N^v_u \subset \mathbf{M}^v_u$ of local minima of $f$.

\begin{prop}
  Let $u \in H^1(\Sigma; \Z/2) - \{0\}$, let $v \in H^2(\Sigma; \Z/2)
  = \Z/2$ and let $P^+$ and $P^-$ be the Abelian varieties associated
  to the double cover of $\Sigma$, given by $u$ as above.  Then, the
  subvariety $N^v_u \subset \mathbf{M}^v_u$ of local minima of $f$ is
  equal to $P^+$ and $P^-$, respectively, for the two values of $v$.
  \qed
\end{prop}

Consequently, $N^v_u$ is connected and, from the properness of $f$, it
follows that $\mathbf{M}^v_u$ is connected, finally finishing the
proof of \thmref{thm:components}.

\providecommand{\bysame}{\leavevmode\hbox to3em{\hrulefill}\thinspace}

\noindent
\textit{\nocorr
  Departamento de Matem{\'a}tica Pura\\
  Faculdade de Ci{\^e}ncias da Universidade do Porto \\
  4099-002 Porto \\
  Portugal \\
  E-mail:} \texttt{pbgothen@fc.up.pt}

\end{document}